\documentclass[]{interact}
\usepackage{enumerate,varwidth,caption,float,color}
\usepackage{graphicx} %use graph format
\usepackage{epstopdf}%To incorporate .eps illustrations using PDFLaTeX, etc.
\usepackage[numbers,sort&compress]{natbib}% Citation support using natbib.sty
\bibpunct[, ]{[}{]}{,}{n}{,}{,}% Citation support using natbib.sty
\makeatletter% @ becomes a letter
\def\NAT@def@citea{\def\@citea{\NAT@separator}}% Suppress spaces between citations using natbib.sty
\makeatother% @ becomes a symbol again

\theoremstyle{plain}% Theorem-like structures provided by amsthm.sty
\newtheorem{theorem}{Theorem}[section]

\newtheorem{proposition}{Proposition}[section]%[theorem]

\theoremstyle{definition}

\theoremstyle{remark}
\newtheorem{remark}{Remark}

\begin{document}
\articletype{ARTICLE TEMPLATE}% Specify the article type or omit as appropriate
\title{Bivariate first-order random coefficient integer-valued autoregressive processes based on modified negative binomial operator}
%Taylor \& Francis \LaTeX\ template for authors (\textsf{Interact} layout + NLM reference style)
\author{
\name{Yixuan Fan\textsuperscript{a}, Dehui Wang\textsuperscript{b*}\thanks{*CONTACT Dehui Wang. Email: wangdehui@lnu.edu.cn}}
\affil{
\textsuperscript{a}School of Mathematics, Jilin University, Changchun, 130012, China;\\
\textsuperscript{b}School of Mathematics and Statistics, Liaoning University, Shenyang, 110036, China}
}
\maketitle

\begin{abstract}
In this paper, a new bivariate random coefficient integer-valued autoregressive process based on modified negative binomial operator with dependent innovations is proposed. Basic probabilistic and statistical properties of this model are derived. To estimate unknown parameters, Yule-Walker, conditional least squares and conditional maximum likelihood methods are considered and evaluated by Monte Carlo simulations. Asymptotic properties of the estimators are derived. Moreover, coherent forecasting and possible extension of the proposed model is provided. Finally, the proposed model is applied to the monthly crime datasets and compared with other models.
\end{abstract}

\section{Introduction}
Integer-valued time series of counts are frequently encountered in real-life situations, such as seismology, certain diseases in epidemiology, criminal incidents and traffic accidents. Al-Osh and Alzaid (1987) proposed the classic first-order integer-valued autoregressive (INAR(1)) process with Poisson marginal to fit time series of counts. To make the INAR(1) model more flexible, Zheng et al. (2007) pointed out that, the autoregressive parameter in the INAR(1) model should vary with time randomly and thus they proposed the random coefficient INAR(1) model. Ding and Wang (2016) incorporated the explanatory variables into the autoregressive coefficient using logit transformation for the INAR(1) process. Yang et al. (2021) developed a self-exciting integer-valued threshold autoregressive process with the random coefficient driven by explanatory variables. To capture the driving effect of covariates on the finite-range time series of counts, Li et al. (2023) introduced a $p$th-order random coefficient mixed binomial autoregressive process with explanatory variables. In the above references, the Bernoulli counting series in the INAR(1) process take only two possible values. The expectation and the variance of the Poisson distributed innovations are equal and they are not always suitable for modelling overdispersed data. Therefore, Risti\'{c} et al. (2009) defined the negative binomial operator for the INAR(1) process with geometric marginal distribution. Afterthat, INAR models based on negative binomial operator attract more attention in the fields of statistics and economics.

To extend the INAR(1) model to the multivariate situation and handle the serial correlation of the successive count data for time series models, Pedeli and Karlis (2011) first proposed the bivariate INAR(1) (BINAR(1)) model,
\begin{align*}
\textit{\textbf{X}}_{t}=\textit{\textbf{A}}\circ\textit{\textbf{X}}_{t-1}+\textit{\textbf{Z}}_{t}=\left[\begin{matrix}\alpha_{1}&0\\0&\alpha_{2}\\\end{matrix}\right]
\circ
\left[\begin{matrix}X_{1,t-1}\\X_{2,t-1}\\\end{matrix}\right]+\left[\begin{matrix}Z_{1,t}\\Z_{2,t}\\
\end{matrix}\right],\; t\in\mathbb{N},
\end{align*}
where $``\circ"$ is the matricial thinning operator and $\textit{\textbf{A}}\circ$ acts as the usual matrix multiplication and keeps the properties of the binomial thinning operator. The innovation $\textit{\textbf{Z}}_{t}$ is a sequence of independent and identically distributed bivariate integer-valued random vectors, which is independent of $\textit{\textbf{X}}_{s}$ for $s<t$. Subsequently, Risti\'{c} et al. (2012) proposed a bivariate integer-valued time series model with geometric marginals based on the negative binomial operator. Pedeli and Karlis (2013a) and Pedeli and Karlis (2013b) considered the parameter estimation problems for BINAR(1) models with bivariate Poisson innovations. Karlis and Pedeli (2013) considered a BINAR(1) process where cross-correlation is constructed by a copula function for the bivariate negative binomial distributed innovations. For modelling count data with bounded range, Scotto et al. (2014) proposed a bivariate binomial autoregressive model which allows for positive and negative cross-correlations marginally. Popovi\'{c} et al. (2016) constructed a BINAR(1) model using the bivariate geometric innovations with different parameters. Popovi\'{c} (2015) and Yu et al. (2020) introduced two different types of random coefficient BINAR(1) processes based on binomial thinning operator with dependent innovations. Su and Zhu (2021) incorporated explanatory variables into BINAR(1) models with different bivariate negative binomial innovations. To better modelling the bivariate and multivariate count time series showing piecewise phenomena, Yang et al. (2023a) proposed a bivariate threshold Poisson INAR(1) process. Then Yang et al. (2023b) introduced a class of multivariate threshold INAR(1) processes driven by explanatory variables to capture the correlations for the counts of time series.

To introduce new members for the BINAR(1) models, Zhang et al. (2020) proposed the extended negative binomial operator $``\diamond"$, written as
\begin{align}\label{eq1}
\alpha\diamond X=\sum_{j=1}^{X+1}G_j,
\end{align}
where $\{G_j\}$ is a sequence of independent and identically distributed geometric random variables with parameter $\frac{\alpha}{1+\alpha}$. It satisfies $\alpha\diamond 0\neq 0$ for all values of $\alpha$. Thereafter, Aleksi\'{c} and Risti\'{c} (2021) renamed this operator as modified negative binomial operator and constructed a new minification INAR(1) model based on it. Furthermore, Qian and Zhu (2022) proposed a new min-INAR(1) process driven by the explanatory variables, which makes the model more practical. However, these present studies draw primarily on the work of fixed operator. In fact, the autoregressive coefficient $\alpha$ in (\ref{eq1}) could vary with time randomly. For instance, let $X_t=\alpha\diamond X_{t-1}+\epsilon_t$ represents the number of crimes under a given area in time $t$, where $\alpha\diamond X_{t-1}$ stands for the occurrence of crimes in time $t-1$ and $\epsilon_t$ denotes the arrival of newly crimes in time $t$. The fixed operator cannot introduce the influence from politics, economics, and other factors which may have a great influence on the observation $X_t$. Therefore, it may be more reasonable to replace the fixed coefficient $\alpha$ with the random coefficient $\alpha_t$. More specifically, the motivation of this model comes from the Beta-negative binomial distribution. Generally speaking, integer-valued time series models can be sorted into two categories: “thinning” operator based models and auto-regression based models. A large fraction of existing works in the literature have focused on modelling continuous-valued time series with certain heavy-tailed distributions. However, current literature on modelling discrete time series with heavy-tailedness is less considered. To model integer-valued time series data with heavy-tailedness, Qian et al. (2020) considered the INAR(1) process based on binomial thinning operator with generalized Poisson-inverse Gaussian innovations. Gorgi (2020) proposed a general class of Beta-negative binomial auto-regression models. Based on the above references, we proposed this new random coefficient BINAR(1) process to further enhance model flexibility.

%To make the BINAR(1) process more applicable, we introduce a new random coefficient BINAR(1) process via replacing $\alpha_i$ with random variable $\alpha_{i,t}\ i=1,2$, and discuss its probabilistic and statistical properties. 
%The motivation for this paper is to introduce a new random coefficient BINAR(1) process based on the modified negative binomial thinning operator 

%2013Flexible Bivariate INAR(1) Processes Using Copulas: The main focus in this article is on models with bivariate Poisson and bivariate negative binomial innovations appropriate for equidispersed and overdispersed data, respectively.

%forecasting=prediction related
The outline of the paper is as follows. In Section \ref{SEC2}, we propose a new bivariate random coefficient INAR(1) process based on modified negative binomial operator and discuss its basic probabilistic and statistical properties. In Section \ref{SEC3}, we consider the Yule-Walker, conditional least squares and conditional maximum likelihood methods to estimate the parameters. In Section \ref{SEC4}, we evaluate the finite-sample performance of these three estimation methods by Monte Carlo simulations respectively. In Section \ref{SEC5}, we discuss the coherent prediction of our proposed random coefficient BINAR(1) model. In Section \ref{SEC6}, we apply the proposed model to a real example to demonstrate its usefulness and flexibility. In Section \ref{SEC7}, we provide the possible extension for our proposed model and Section \ref{SEC8} concludes.
%investigate

\section{Random coefficient BINAR(1) processes based on modified negative binomial operator}\label{SEC2}
In this section, we propose a new bivariate random coefficient INAR(1) process, called BRCMNBINAR(1). The process is defined as
\begin{align}\label{eq2}
\textit{\textbf{X}}_{t}=\textit{\textbf{A}}_t\diamond\textit{\textbf{X}}_{t-1}+\textit{\textbf{Z}}_{t}=\left[\begin{matrix}\alpha_{1,t}&0\\0&\alpha_{2,t}\\
\end{matrix}\right]
\diamond
\left[\begin{matrix}X_{1,t-1}\\X_{2,t-1}\\\end{matrix}\right]+\left[\begin{matrix}Z_{1,t}\\Z_{2,t}\\
\end{matrix}\right],\; t\in\mathbb{N},
\end{align}
where 
\begin{enumerate}[(i)]%which take values on $(0,1)$
\item For $i=1,2$, we assume that $\{\alpha_{i,t}\}$ is a sequence of independent and identically distributed random variables, follows the Beta prime distribution with parameters $m_i$ and known constant $l$. The probability density function of $\alpha_{i,t}$ is given by $f(x)=x^{m_i-1}(1+x)^{-m_i-l-1}/B(m_i,l+1)$, where $B(\cdot)$ denotes the Beta function. $\{\alpha_{1,t}\}$ and $\{\alpha_{2,t}\}$ are mutually independent. Assume that $l > 1$ and $m_i < l$, the expectation and variance of $\alpha_{i,t}$ are given by
\begin{align*}
E(\alpha_{i,t})=\alpha_i=\frac{m_i}{l},\ Var(\alpha_{i,t})=\sigma_{\alpha_i}^2=\frac{m_i(m_i+l)}{(l-1)l^2},\ i=1,2.
\end{align*}

\item The modified negative binomial operator $``\diamond"$ is defined in (1).
%proposed by Zhang et al. (2020), is defined as $\alpha_{i,t}\diamond X_{t-1}=\sum_{j=1}^{X_{t-1}+1}G_{ij}^{(t)}$, where $\alpha_{i,t}\in (0,1)$, $\{G_{ij}^{(t)}\}$ is a sequence of independent and identically distributed geometric random variables with parameter $\frac{\alpha_{i,t}}{1+\alpha_{i,t}}$, $i=1,2$.

\item $\textit{\textbf{A}}_t=\left[\begin{matrix}\alpha_{1,t}&0\\0&\alpha_{2,t}\\\end{matrix}\right]$ is a random matrix and $\textit{\textbf{A}}_t\diamond$\ is the random matricial operation which acts as the usual matrix multiplication and keeps the properties of the random coefficient operator. %For $i=1,2$, $\forall t\in\mathbb{N}$, $\alpha_{1,t}\diamond X_{1,t-1}$ and $\alpha_{2,t}\diamond X_{2,t-1}$ are mutually independent.

\item $\{\textit{\textbf{Z}}_{t}\}_{t\in\mathbb{N}}$ is a sequence of independent and identically distributed bivariate non-negative integer-valued random vectors, with joint probability mass function $f_\textbf{Z}(x,y)$. For fixed $t$, $\textit{\textbf{Z}}_{t}$ is independent of $\textit{\textbf{A}}_t\diamond\textit{\textbf{X}}_{t-1}$ and $\textit{\textbf{X}}_{s}$ for $s<t$. We assume that the expectation and variance of $Z_{i,t}$ are $E(Z_{i,t})=\mu_{i},\ Var(Z_{i,t})=\sigma^2_{{Z_i}},\ i=1,2$. The dependence between the innovations $Z_{1,t}$ and $Z_{2,t}$ is described by $Cov(Z_{1,t},Z_{2,t})=\phi$. 
\end{enumerate}
\begin{remark}
Similar to Li et al. (2018), the value of $l$ can be selected by the decision maker based on experience. For example, assume that $1<\underline{l}\leqslant l\leqslant\bar{l}$, the value of $l$ can be searched over in $[\underline{l},\bar{l}]$ by maximizing the likelihood function. %We can also use the two-step CLS method proposed in Hwang and Basawa (1998) to calculate $\hat{l}$, further details are discussed in Section \ref{SEC3}.
\end{remark}

Model (2) is a generalization of the model proposed in Zhang et al. (2020). It can be seen as a special case of our proposed model when $\sigma_{\alpha_i}^2=0,\ i=1,2$. Notice that the $i$th element of model (\ref{eq2}) can be written as $X_{i,t}=\alpha_{i,t}\diamond X_{i,t-1}+Z_{i,t},\ i=1,2$. Furthermore, the dependence between these two series is introduced via the dependence between the innovations $\{Z_{1,t}\}$ and $\{Z_{2,t}\}$. Throughout the paper, we assume that $(m_1,\mu_1)\neq (m_2,\mu_2)$.

In the following proposition, we provide the transition probabilities of the BRCMNBINAR(1) process. 
\begin{proposition}\label{proposition2.1}
The process $\{\textbf{\textit{X}}_t\}_{t\in\mathbb{N}}$ defined in (\ref{eq2}) is a Markov process on $\mathbb{N}_0^2$ with the following transition probabilities:   
\begin{align}\label{eq3}
&P(\textbf{\textit{X}}_{t}=\textbf{\textit{x}}_{t}|\textbf{\textit{X}}_{t-1}=\textbf{\textit{x}}_{t-1})\nonumber\\
=&P(X_{1,t}=x_{1,t},X_{2,t}=x_{2,t}|X_{1,t-1}=x_{1,t-1},X_{2,t-1}=x_{2,t-1})\nonumber\\
=&\sum_{k=0}^{x_{1,t}}\sum_{s=0}^{x_{2,t}}f_\textbf{Z}(x_{1,t}-k,x_{2,t}-s)\binom{x_{1,t-1}+k}{k}\frac{\Gamma(x_{1,t-1}+l+2)\Gamma(k+m_1)\Gamma(l+m_1+1)}{\Gamma(x_{1,t-1}+k+m_1+l+2)\Gamma(l+1)\Gamma(m_1)}\nonumber\\
&\qquad\qquad\times\binom{x_{2,t-1}+s}{s}\frac{\Gamma(x_{2,t-1}+l+2)\Gamma(s+m_2)\Gamma(l+m_2+1)}{\Gamma(x_{2,t-1}+s+m_2+l+2)\Gamma(l+1)\Gamma(m_2)}.
\end{align}
\end{proposition}

For the random matrix $\textbf{\textit{A}}_t$, we have $E(\textbf{\textit{A}}_t)=\textbf{\textit{A}}=\left[\begin{matrix}\alpha_{1}&0\\0&\alpha_{2}\\\end{matrix}\right]$. The following proposition states some properties for the matrix $\textbf{\textit{A}}_t$.
\begin{proposition}\label{proposition2.2}
Let $\textbf{\textit{A}}_t\diamond\textbf{\textit{X}}_{t-1}$ be a random matricial operation, $\textbf{\textit{Y}}$ is a bivariate random vector which is independent of the random matricial operation, then we have
\begin{enumerate}[(i)]
\item $E(\textbf{\textit{A}}_t\diamond\textbf{\textit{X}}_{t-1})=\textbf{\textit{A}}E(\textbf{\textit{X}}_{t-1}+\textbf{\textit{e}})$,

\item $E[(\textbf{\textit{A}}_t\diamond\textbf{\textit{X}}_{t-1})\textbf{Y}^T]=\textbf{\textit{A}}[E(\textbf{\textit{X}}_{t-1})+\textbf{\textit{e}}]E(\textbf{\textit{Y}}^T)$,

\item $E[\textbf{Y}(\textbf{\textit{A}}_t\diamond\textbf{X}_{t-1})^T]=E(\textbf{\textit{Y}})[E(\textbf{\textit{X}}_{t-1}^T)+\textbf{\textit{e}}^T]\textbf{\textit{A}}$,

\item $E[(\textbf{\textit{A}}_t\diamond\textbf{X}_{t-1})(\textbf{\textit{A}}_t\diamond\textbf{X}_{t-1})^T]=\textbf{\textit{A}}E[(\textbf{\textit{X}}_{t-1}+\textbf{\textit{e}})(\textbf{\textit{X}}_{t-1}+\textbf{\textit{e}})^T]\textbf{\textit{A}}+\textbf{\textit{C}}$,
\end{enumerate}
where $\textbf{\textit{e}}=(1,1)^T$, $c_{ii}=(\sigma^2_{\alpha_i}+\alpha_i^2+\alpha_i)(EX_{i,t-1}+1)+\sigma^2_{\alpha_i}E(X_{i,t-1}+1)^2,\ i=1,2$ and $c_{12}=c_{21}=m_1m_2\phi/(l^2-m_1m_2)$.
\iffalse
\begin{proof}\rm
These properties are easied to varified and we omit the details.
\end{proof}
\fi
\end{proposition}

%先给平稳性
In the following proposition, the strict stationarity and ergodicity for the process (\ref{eq2}) are established, which are useful to obtain the asymptotic properties of the parameter estimators. 
\begin{proposition}\label{proposition2.3}
If $0<\alpha_i^2+\sigma_{\alpha_i}^2<1,\ i=1,2$, then there exists a unique strictly stationary and ergodic process $\{\textbf{\textit{X}}_t\}_{t\in\mathbb{N}}$ satisfying (\ref{eq2}). 
\end{proposition}

The moments and conditional moments will be useful in parameter estimation in Section \ref{SEC3}. For BRCMNBINAR(1) process, we have the following results.
\begin{proposition}\label{proposition2.4}
Let $\{\textbf{\textit{X}}_t\}_{t\in\mathbb{N}}$ be a strictly stationary process given by (\ref{eq2}), then for $t\geqslant 1,\ k\geqslant 0,\ i,j=1,2$ and $i\neq j$, we have
\begin{enumerate}[(i)]
\item $E(X_{i,t+1}|X_{1,t},X_{2,t})=\alpha_i(X_{i,t}+1)+\mu_{i}$,\\
$Var(X_{i,t+1}|X_{1,t},X_{2,t})=\sigma^2_{\alpha_i}(X_{i,t}+1)^2+(\sigma^2_{\alpha_i}+\alpha_i^2+\alpha_i)(X_{i,t}+1)+\sigma^2_{Z_i}$,\\ 
$Cov(X_{1,t+1},X_{2,t+1}|X_{1,t},X_{2,t})=\phi$,%check

\item $E(X_{i,t})=(m_i+l\mu_{i})/(l-m_i)$,\\
$Var(X_{i,t})=\big[(\sigma^2_{\alpha_i}+\alpha_i^2+\alpha_i)\frac{l(\mu_{i}+1)}{(l-m_i)}+\frac{m_i(m_i+l)(\mu_{i}+1)^2}{(l-1)(l-m_i)^2}+\sigma^2_{Z_i}\big]/(1-\alpha_i^2-\sigma^2_{\alpha_i})
$,\\
$Cov(X_{i,t+k},X_{i,t})=\alpha_i^kVar(X_{i,t})$,\\
$Corr(X_{i,t+k},X_{i,t})=\alpha_i^k$,%check

\item $Cov(X_{1,t},X_{2,t})=\frac{l^2\phi}{l^2-m_1m_2}$,\\
$Cov(X_{i,t+k},X_{j,t})=m_i^k\phi/(l^{k}-l^{k-2}m_1m_2)$.
\end{enumerate}
\end{proposition}

\section{Parameter estimation}\label{SEC3}
In this section, we consider parameter estimation problem using Yule-Walker (YW), conditional least squares (CLS) and conditional maximum likelihood (CML) methods. Suppose that $\textbf{\textit{X}}_1,\textbf{\textit{X}}_2,\cdots,\textbf{\textit{X}}_n$ is a sequence of observations generated by the BRCMNBINAR(1) process. We consider parameter estimation of $\pmb\theta=(m_1,m_2,\mu_1,\mu_2,\phi)^T$ and the corresponding variance $\pmb{\eta}=(\sigma^2_{\alpha_1},\sigma^2_{\alpha_2},\sigma^2_{Z_1},\sigma^2_{Z_2})^T$ respectively. 
\subsection{Yule-Walker estimation}
Let $\gamma_i(0)=Var(X_{i,t}),\ \gamma_i(1)=Cov(X_{i,t},X_{i,t-1}),\ \gamma_{ij}(0)=Cov(X_{1,t},X_{2,t})$. Then we have 
\begin{align*}
&\hat{\gamma}_i(0)=\frac{1}{n}\sum_{t=1}^{n}(X_{i,t}-\frac{1}{n}\sum_{t=1}^{n}X_{i,t})^2,\\   
&\hat{\gamma}_i(1)=\frac{1}{n-1}\sum_{t=2}^{n}(X_{i,t}-\frac{1}{n}\sum_{t=1}^{n}X_{i,t})(X_{i,t-1}-\frac{1}{n}\sum_{t=1}^{n}X_{i,t}),\\
&\hat{\gamma}_{ij}(0)=\frac{1}{n}\sum_{t=1}^{n}(X_{1,t}-\frac{1}{n}\sum_{t=1}^{n}X_{1,t})(X_{2,t}-\frac{1}{n}\sum_{t=1}^{n}X_{2,t}).
\end{align*}
Thus, we have the following YW-estimator 
\begin{align*}
&\hat{m}_i^{YW}=\hat{\gamma}_i(1)l/\gamma_i(0),\ \hat{\mu}_i^{YW}=\frac{1}{n-1}\sum_{t=2}^{n}X_{i,t}-\hat{m}_i^{YW}(X_{i,t-1}+1)/l,\ i=1,2,\\
&\hat{\phi}^{YW}=(1-\hat{m}_1^{YW}\hat{m}_2^{YW}/l^2)\hat{\gamma}_{ij}(0).
\end{align*}

%=(\hat{m}_1^{YW},\hat{m}_2^{YW},\hat{\mu}_1^{YW},\hat{\mu}_2^{YW},\hat{\phi}^{YW})'
%easily 
The strong consistency of the YW-estimator $\hat{\pmb\theta}^{YW}$ can be proved along the similar way of Theorem 4.1 in Liu et al. (2015) and Theorem 1 in Yu et al. (2020), thus we omit the proof here. %In addition, its asymptotic distribution is equivalent to the CLS-estimator $\hat{\pmb\theta}^{CLS}$ by taking analogous argument in Pedeli and Karlis (2013a).

\subsection{Conditional least squares estimation}
Similar to Klimko and Nelson (1978), we first consider the conditional least squares estimator of the parameters $\pmb{\gamma}=(m_1,m_2,\mu_1,\mu_2)^T$. Let 
\begin{align*}
g(\pmb{\gamma},\textit{\textbf{X}}_{t-1})=E(\textit{\textbf{X}}_t|\textit{\textbf{X}}_{t-1})=\textit{\textbf{A}}(\textit{\textbf{X}}_{t-1}+\textit{\textbf{e}})+E(\textit{\textbf{Z}}_{t}),\ \textit{\textbf{u}}_t=\textit{\textbf{X}}_{t}-g(\pmb{\gamma},\textit{\textbf{X}}_{t-1}),
\end{align*}
where $E(\textit{\textbf{Z}}_{t})=(\mu_1,\mu_2)^T$. Then the CLS-estimator $\hat{\pmb{\gamma}}^{CLS}$ is obtained by minimizing the following function
\begin{align*}
S_1(\pmb{\gamma})
&=\sum_{t=2}^{n}\textbf{\textit{u}}_t^T\textbf{\textit{u}}_t\\
&=\sum_{t=2}^{n}(X_{1,t}-m_1(X_{1,t-1}+1)/l-\mu_1)^2+(X_{2,t}-m_2(X_{2,t-1}+1)/l-\mu_2)^2.
\end{align*}

The solution of equation $D_t(\pmb{\gamma})=-(1/2)\partial S_1(\pmb{\gamma})/\partial{\pmb{\gamma}}=\textbf{\textit{0}}$ yields the CLS-estimator of $\pmb{\gamma}$. It is easy to verify that $D_t(\pmb{\gamma})=(D_{t1}(\pmb{\gamma}),D_{t2}(\pmb{\gamma}),D_{t3}(\pmb{\gamma}),D_{t4}(\pmb{\gamma}))^T$ with 
\begin{align*}
&D_{t1}(\pmb{\gamma})=-\frac{1}{2}\frac{\partial S_1(\pmb{\gamma})}{\partial{m_1}}=\frac{1}{l}\sum_{t=2}^{n}(X_{1,t}-m_1(X_{1,t-1}+1)/l-\mu_1)(X_{1,t-1}+1),\\
&D_{t2}(\pmb{\gamma})=-\frac{1}{2}\frac{\partial S_1(\pmb{\gamma})}{\partial{m_2}}=\frac{1}{l}\sum_{t=2}^{n}(X_{2,t}-m_2(X_{2,t-1}+1)/l-\mu_2)(X_{2,t-1}+1),\\
&D_{t3}(\pmb{\gamma})=-\frac{1}{2}\frac{\partial S_1(\pmb{\gamma})}{\partial{\mu_1}}=\sum_{t=2}^{n}(X_{1,t}-m_1(X_{1,t-1}+1)/l-\mu_1),\\
&D_{t4}(\pmb{\gamma})=-\frac{1}{2}\frac{\partial S_1(\pmb{\gamma})}{\partial{\mu_2}}=\sum_{t=2}^{n}(X_{2,t}-m_2(X_{2,t-1}+1)/l-\mu_2).
\end{align*}
Then we can obtain the CLS-estimator $\hat{\pmb{\gamma}}^{CLS}$ as follows:
\begin{align*}
\hat{\pmb{\gamma}}^{CLS}=\textit{\textbf{G}}_{n-1}^{-1}\textit{\textbf{b}}_{n-1},
\end{align*}
where
\begin{align*}
\textit{\textbf{G}}_{n-1}=   
\begin{pmatrix}
\frac{1}{l}\sum_{t=2}^{n}(X_{1,t-1}+1)^2&0&\sum_{t=2}^{n}(X_{1,t-1}+1)&0,\\
0&\frac{1}{l}\sum_{t=2}^{n}(X_{2,t-1}+1)^2&0&\sum_{t=2}^{n}(X_{2,t-1}+1)\\
\frac{1}{l}\sum_{t=2}^{n}(X_{1,t-1}+1)&0&n-1&0\\
0&\frac{1}{l}\sum_{t=2}^{n}(X_{2,t-1}+1)&0&n-1
\end{pmatrix},
\end{align*}  
and 
\begin{align*}
\textit{\textbf{b}}_{n-1}=(\sum_{t=2}^{n}X_{1,t}(X_{1,t-1}+1),\sum_{t=2}^{n}X_{2,t}(X_{2,t-1}+1),\sum_{t=2}^{n}X_{1,t},\sum_{t=2}^{n}X_{2,t})^T.
\end{align*}

The following result states the strong consistency and asymptotic normality of the CLS-estimator $\hat{\pmb{\gamma}}^{CLS}$.
\begin{theorem}\label{th3.1}
Let $\{\textit{\textbf{X}}_t\}_{t\in\mathbb{N}}$ be the BRCMNBINAR(1) process defined in (\ref{eq2}). Assume that $E|X_{i,t}|^4<\infty,\ i=1,2$, the CLS-estimator $\hat{\pmb{\gamma}}^{CLS}=(\hat{m}^{CLS}_1,\hat{m}^{CLS}_2,\hat{\mu}^{CLS}_1,\hat{\mu}^{CLS}_2)^T$ is strongly consistent and asymptotically normal, that is
\begin{align}
\sqrt{n-1}(\hat{\pmb{\gamma}}^{CLS}-\pmb{\gamma})\xrightarrow{L}N(\textit{\textbf{0}},{\textit{\textbf{V}}_1}^{-1}{\textit{\textbf{W}}_1}{\textit{\textbf{V}}_1}^{-1}),
\end{align} 
where $\textit{\textbf{V}}_1=E\left(\frac{\partial g(\pmb{\gamma},\textbf{\textit{X}}_{t-1})^T\partial g(\pmb{\gamma},\textbf{\textit{X}}_{t-1})}{\partial\pmb{\gamma}\partial\pmb{\gamma}^T}\right)$ 
and $\textit{\textbf{W}}_1=E\left(\frac{\partial g(\pmb{\gamma},\textbf{\textit{X}}_{t-1})^T\textbf{\textit{u}}_t\textbf{\textit{u}}_t^T\partial g(\pmb{\gamma},\textbf{\textit{X}}_{t-1})}{\partial{\pmb{\gamma}}\partial\pmb{\gamma}^T}\right)$, $\xrightarrow{L}$ denotes convergence in distribution.
\end{theorem}

We now consider the problem of estimating $\phi$. Let $Y_t=(X_{1,t}-E(X_{1,t}|\textit{\textbf{X}}_{t-1}))(X_{2,t}-E(X_{2,t}|\textit{\textbf{X}}_{t-1}))$. From proposition \ref{proposition2.4}, we have $E(Y_t|\textit{\textbf{X}}_{t-1})=Cov(X_{1,t},X_{2,t}|X_{1,t-1},X_{2,t-1})=\phi$. Then we construct the following criterion function
\begin{align*}
S_2(\phi)&=\sum_{t=2}^{n}(Y_t-E(Y_t|\textit{\textbf{X}}_{t-1}))^2\\
&=\sum_{t=2}^{n}[(X_{1,t}-E(X_{1,t}|\textit{\textbf{X}}_{t-1}))(X_{2,t}-E(X_{2,t}|\textit{\textbf{X}}_{t-1}))-\phi]^2.
%&=\sum_{t=2}^{n}[(X_{1,t}-\hat{m}^{CLS}_1(X_{1,t-1}+1)/l-\hat{\mu}^{CLS}_1)(X_{2,t}-\hat{m}^{CLS}_2(X_{2,t-1}+1)/l-\hat{\mu}^{CLS}_2)-\phi]^2.
\end{align*}

\begin{theorem}\label{th3.2}
If $E|X_{i,t}|^4<\infty,\ i=1,2$, then the CLS-estimator $\hat{\phi}^{CLS}$ is strongly consistent and asymptotically normal, that is
\begin{align*}
\sqrt{n-1}(\hat{\phi}^{CLS}-\phi)\xrightarrow{L} N(0,\delta^2),
\end{align*}  
where $\delta^2=E[(X_{1,t}-m_1(X_{1,t-1}+1)/l-\mu_1)(X_{2,t}-m_2(X_{2,t-1}+1)/l-\mu_2)-\phi]^2$.
\end{theorem}

Now we consider the problem of estimating the variance parameter $\pmb{\eta}=(\sigma^2_{\alpha_1},\sigma^2_{\alpha_2},\sigma^2_{Z_1},\sigma^2_{Z_2})^T$. We use the two-step conditional least squares method introduced by Karlsen and Tj\o{}stheim (1988) to estimate the parameter $\pmb\eta$. Let 
\begin{align}\label{eq11}
S_3(\pmb{\eta})=&\sum_{t=2}^{n}\textbf{\textit{U}}_t^T\textbf{\textit{U}}_t\nonumber\\
=&\sum_{t=2}^{n}\sum_{j=1}^{2}[(X_{j,t}-E(X_{j,t}|X_{j,t-1}))^2-Var(X_{j,t}|X_{j,t-1})]^2\nonumber\\
=&\sum_{t=2}^{n}\sum_{j=1}^{2}[(X_{j,t}-\hat{\alpha}_j(X_{j,t-1}+1)-\hat{\mu}_j)^2-\sigma^2_{\alpha_j}(X_{j,t}+1)^2\nonumber\\
&-(\sigma^2_{\alpha_j}+\hat{\alpha}_j^2+\hat{\alpha}_j)(X_{j,t}+1)-\sigma^2_{Z_j}]^2.
\end{align}

Similar to Hwang and Basawa (1998), we minimize (\ref{eq11}) with $\hat{\alpha}_i$ replaced by $\hat{m}_i^{CLS}/l$ in $\hat{\pmb{\gamma}}^{CLS}$,\ $i=1,2$. Then the CLS-estimator $\hat{\pmb{\eta}}^{CLS}$ can be numerically calculated by minimizing the criterion function:
\begin{align*}
\hat{\pmb{\eta}}^{CLS}=\arg\ \min_{\pmb{\eta}}S_3(\pmb{\eta}).
\end{align*}

Solving the equation system $-(1/2)\partial S_3(\pmb\eta)/\partial{\pmb\eta}=\textbf{\textit{0}}$, we obtain the CLS-estimator of $\pmb\eta$ as follows:
\begin{align*}
\hat{\pmb{\eta}}^{CLS}=\textbf{\textit{H}}_{n-1}^{-1}\textbf{\textit{d}}_{n-1},  
\end{align*}
where 
\begin{align*}
\textbf{\textit{d}}_{n-1}=
\begin{pmatrix}
\sum_{t=2}^{n}[(X_{1,t}-\hat{\alpha}_1(X_{1,t-1}+1)-\hat{\mu}_1)^2-(\hat{\alpha}_1^2+\hat{\alpha}_1)(X_{1,t}+1)][(X_{1,t}+1)^2+X_{1,t}+1]\\
\sum_{t=2}^{n}[(X_{2,t}-\hat{\alpha}_2(X_{2,t-1}+1)-\hat{\mu}_2)^2-(\hat{\alpha}_2^2+\hat{\alpha}_2)(X_{2,t}+1)][(X_{2,t}+1)^2+X_{2,t}+1]\\
\sum_{t=2}^{n}(X_{1,t}-\hat{\alpha}_1(X_{1,t-1}+1)-\hat{\mu}_1)^2-(\hat{\alpha}_1^2+\hat{\alpha}_1)(X_{1,t}+1)\\
\sum_{t=2}^{n}(X_{2,t}-\hat{\alpha}_2(X_{2,t-1}+1)-\hat{\mu}_2)^2-(\hat{\alpha}_2^2+\hat{\alpha}_2)(X_{2,t}+1)
\end{pmatrix},
\end{align*}
and $\textbf{\textit{H}}_{n-1}=$
\begin{align*}
\footnotesize{
\setlength{\arraycolsep}{1.2pt}
\begin{pmatrix}
\sum_{t=2}^{n}[(X_{1,t}+1)^2+X_{1,t}+1]^2&0&\sum_{t=2}^{n}(X_{1,t}+1)^2+X_{1,t}+1&0,\\
0&\sum_{t=2}^{n}[(X_{2,t}+1)^2+X_{2,t}+1]^2&0&\sum_{t=2}^{n}(X_{2,t}+1)^2+X_{2,t}+1\\
\sum_{t=2}^{n}(X_{1,t}+1)^2+X_{1,t}+1&0&n-1&0\\
0&\sum_{t=2}^{n}(X_{2,t}+1)^2+X_{2,t}+1&0&n-1
\end{pmatrix}}.
\end{align*}

\begin{theorem}\label{th3.3}%similar to zhangrui Th3
Let $\{\textit{\textbf{X}}_t\}_{t\in\mathbb{N}}$ be the BRCMNBINAR(1) process defined in (\ref{eq2}). Assume that $E|X_{i,t}|^4<\infty,\ i=1,2$, the CLS-estimator $\hat{\pmb\eta}^{CLS}$ is strongly consistent and asymptotically normal, that is
\begin{align}
\sqrt{n-1}(\hat{\pmb\eta}^{CLS}-\pmb\eta)\xrightarrow{L}N(\textit{\textbf{0}},\textbf{V}_2^{-1}\textbf{W}_2\textbf{V}_2^{-1}),
\end{align}  
where the representations of $\textit{\textbf{V}}_2$ and $\textit{\textbf{W}}_2$ are provided in Appendix.
\end{theorem}

\begin{figure}[!t]
\centering
\includegraphics[width=14cm,height=8cm]{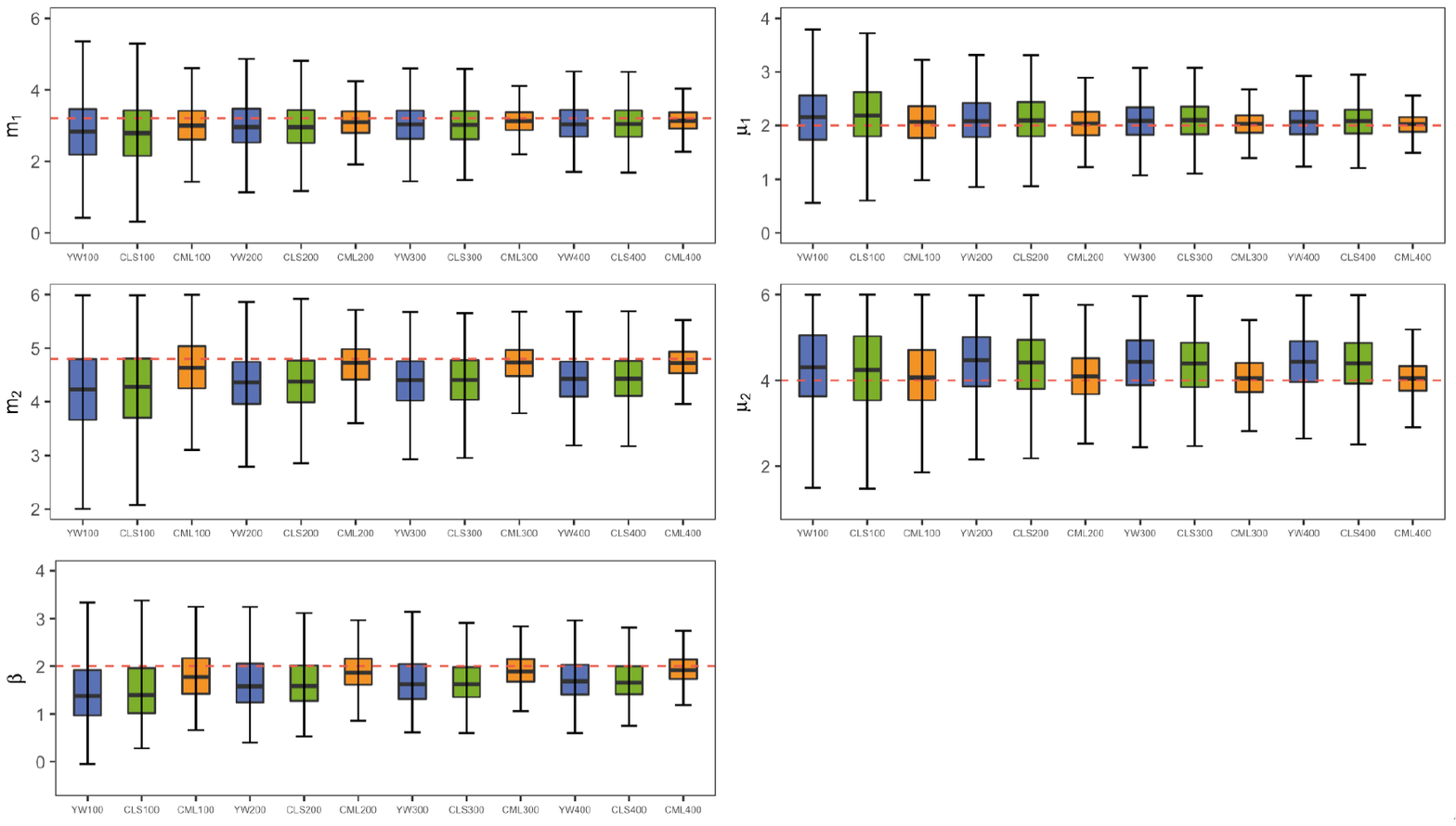}
\caption{The boxplots of the estimates for Scenario A1.}\label{figure1}
\includegraphics[width=14cm,height=8cm]{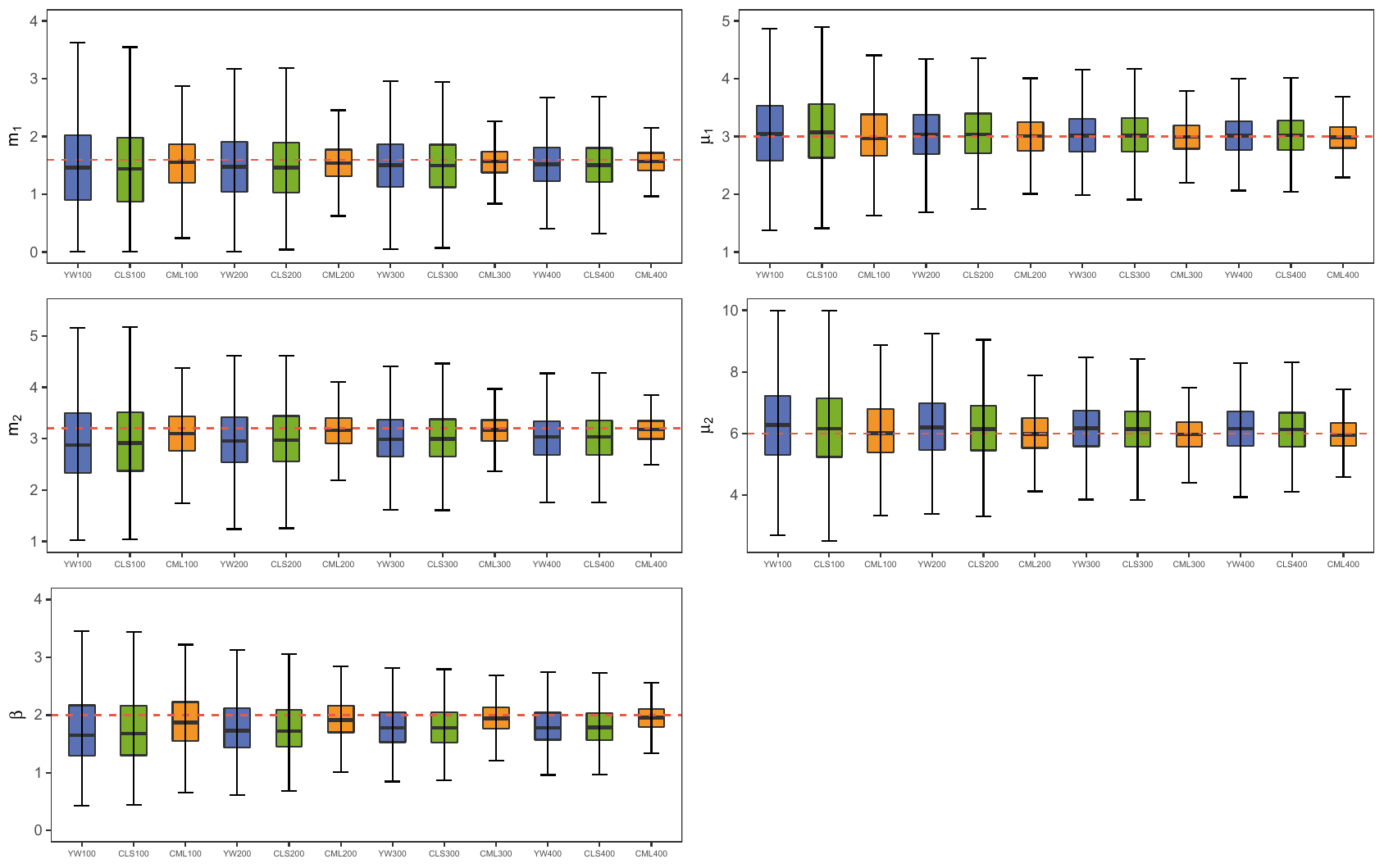}
\caption{The boxplots of the estimates for Scenario A4.}\label{figure2}
\end{figure}

\begin{figure}[!t]
\centering
\includegraphics[width=14cm,height=8cm]{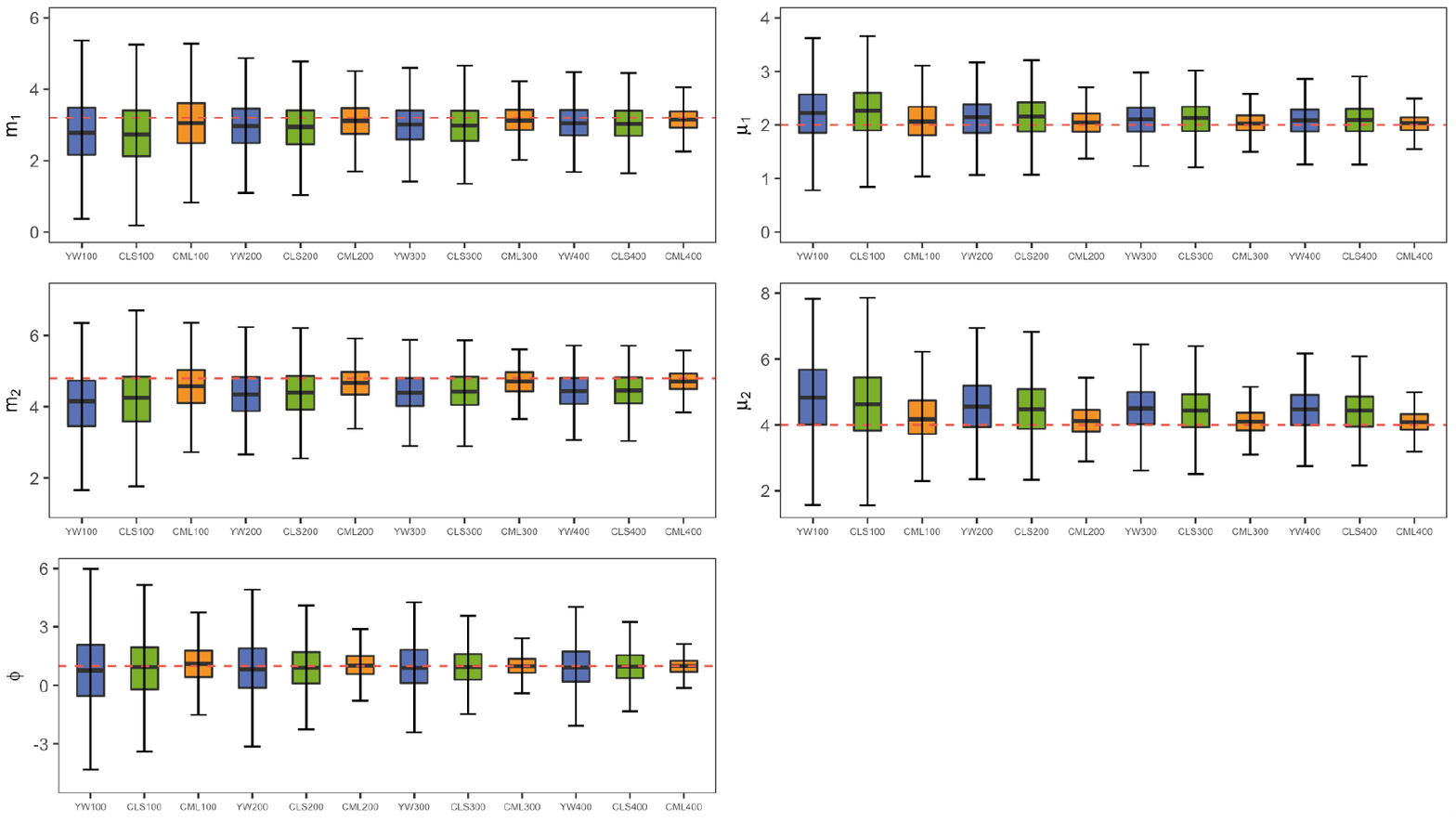}
\caption{The boxplots of the estimates for Scenario B1.}
\label{figure3}
\includegraphics[width=14cm,height=8cm]{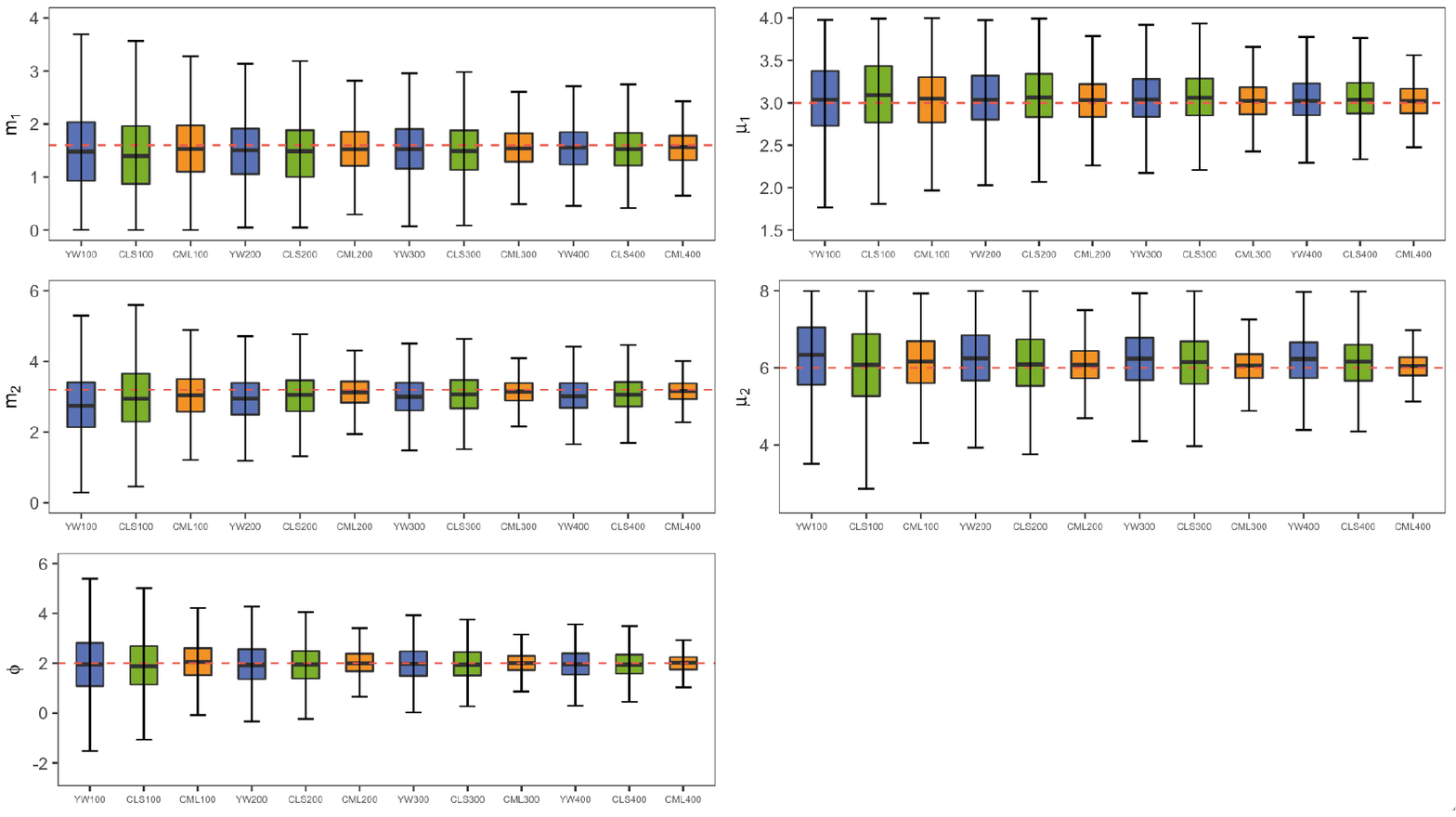}
\caption{The boxplots of the estimates for Scenario B4.}
\label{figure4}
\end{figure}

\subsection{Conditional maximum likelihood estimation}
The conditional likelihood function for the BRCMNBINAR(1) process can be written as 
\begin{align*}
L(\pmb\tau)=\prod_{t=2}^{n}P(\textit{\textbf{X}}_t=\textit{\textbf{x}}_t|\textit{\textbf{X}}_{t-1}=\textit{\textbf{x}}_{t-1}),
\end{align*}
where $\pmb\tau$ is the parameter vector in the transition probabilities $P(\textit{\textbf{X}}_t=\textit{\textbf{x}}_t|\textit{\textbf{X}}_{t-1}=\textit{\textbf{x}}_{t-1})$ given in equation (\ref{eq3}). The CML-estimator $\hat{{\mit\pmb\tau}}^{CML}$ is obtained by maximizing the conditional log-likelihood function: 
\begin{align*}
\hat{{\mit\pmb\tau}}^{CML}=\arg\ \max_{{\pmb\tau}\in\Theta}\ \log L(\pmb\tau).
\end{align*}
Analytical expression for this equation cannot be found. Thus, we need to employ the numerical procedure to solve it. In this paper, we optimize the function by using $\mathrm{nlm}()$ in $\mathcal{R}$ software. The asymptotic normality of the CML-estimator $\hat{{\mit\pmb\tau}}^{CML}$ can be proved by verifying the regularity conditions for estimation of Markov processes in Billingsley (1961). Furthermore, Yang et al. (2023a) provides a formalized proof of the asymptotic behavior of the CML-estimator for the bivariate threshold Poisson integer-valued autoregressive processes.

\section{Simulation studies}\label{SEC4}
In this section, we evaluate the performance of YW-estimator, CLS-estimator and CML-estimator by Monte Carlo simulations. We consider the following bivariate distributions of the innovation $\textit{\textbf{Z}}_t$ for the BRCMNBINAR(1) process.

Case 1: Assume that $(Z_{1,t},Z_{2,t})$ follows the bivariate negative binomial distribution with the joint probability mass function%(BVNB)  
\begin{align*}
P(Z_{1,t}=k,Z_{2,t}=s)=\frac{\Gamma(\beta^{-1}+k+s)\mu_1^k\mu_2^s\beta^{-\beta^{-1}}}{\Gamma(\beta^{-1})\Gamma(k+1)\Gamma(s+1)(\mu_1+\mu_2+\beta^{-1})^{k+s+\beta^{-1}}},
\end{align*}
where $\mu_1,\mu_2,\beta>0$. The bivariate negative binomial distribution is denoted by BVNB($\mu_1,\mu_2,\beta$). If $\beta\rightarrow 0$, then the BVNB distribution will degenerate to the product of two independent Poisson distributions. For short, we refer to this model as I-BRCMNBINAR(1) model.

Case 2: Assume that $(Z_{1,t},Z_{2,t})$ follows the bivariate Poisson distribution with the joint probability mass function 
\begin{align*}
P(Z_{1,t}=k,Z_{2,t}=s)=&e^{-(\mu_1+\mu_2-\phi)}\frac{(\mu_1-\phi)^{k}}{k!}\frac{(\mu_2-\phi)^{s}}{s!}\\
&\times\sum_{i=0}^{min(k,s)}\binom{k}{i}\binom{s}{i}i!(\frac{\phi}{(\mu_1-\phi)(\mu_2-\phi)})^{i},
\end{align*}
where $\mu_1,\mu_2>0,\phi\in[0,\min(\mu_1,\mu_2))$. The bivariate Poisson distribution is denoted by BP($\mu_1,\mu_2,\phi$). The marginal distribution of BP($\mu_1,\mu_2,\phi$) is Poisson distribution with parameters $\mu_1$ and $\mu_2$ respectively. For short, we refer to this model as II-BRCMNBINAR(1) model.
%table1
\begin{table}[!t]
\centering
\caption{Estimated parameters with $n$=100 for scenarios A1-A4.}
\label{table1}
\resizebox{\linewidth}{!}{
\begin{tabular}{rlccccccccccc}
\toprule
&       &  &  YW   &      &       &  &  CLS     &       &       &  &   CML    &  \\
\cmidrule{3-5}\cmidrule{7-9}\cmidrule{11-13}    \multicolumn{1}{l}{Model} & Para. & \multicolumn{1}{l}{bias} & \multicolumn{1}{l}{MSE} & \multicolumn{1}{l}{SD} &       & \multicolumn{1}{l}{bias} & \multicolumn{1}{l}{MSE} & \multicolumn{1}{l}{SD} &       & \multicolumn{1}{l}{bias} & \multicolumn{1}{l}{MSE} & \multicolumn{1}{l}{SD} \\
\midrule
\multicolumn{1}{l}{A1} 
& $m_1=3.2$ & -0.3514 & 0.9667 & 0.9183 &  & -0.3902 & 1.0158 & 0.9293 &       & -0.2002 & 0.4143 & 0.6117 \\
& $m_2=4.8$ & -0.6019 & 1.0942 & 0.8555 &       & -0.5551 & 1.0428 & 0.8571 &  & -0.1619 & 0.3946 & 0.6070 \\
& $\mu_1=2$ & 0.1672 & 0.4007 & 0.6106 &       & 0.2001 & 0.4153 & 0.6126 &       & 0.0817 & 0.1993 & 0.4389 \\
& $\mu_2=4$ & 0.7589 & 2.4135 & 1.3556 &       & 0.6258 & 2.2302 & 1.3560 &       & 0.1924 & 0.8962 & 0.9269 \\
& $\beta=2$ & -0.2747 & 1.3864 & 1.1450 &       & -0.1926 & 6.1753 & 2.4775 &       & -0.1224 & 0.4130 & 0.6309 \\
\cmidrule{1-5}\cmidrule{7-9}\cmidrule{11-13}    \multicolumn{1}{l}{A2} 
& $m_1=3.2$ & -0.8391 & 2.3718 & 1.2914 &       & -0.8458 & 2.3790 & 1.2898 &       & -0.6203 & 1.7652 & 1.1749 \\
& $m_2=4.8$ & -1.3335 & 4.3341 & 1.5987 &       & -1.2744 & 4.2547 & 1.6219 &       & -0.8421 & 3.5668 & 1.6905 \\
& $\mu_1=3$ & -0.2520 & 1.9743 & 1.3823 &       & -0.2522 & 1.9689 & 1.3803 &       & -0.4247 & 1.6341 & 1.2057 \\
& $\mu_2=6$ & -0.0516 & 8.9283 & 2.9876 &       & -0.2658 & 8.5428 & 2.9107 &       & -0.8730 & 6.5121 & 2.3979 \\
& $\beta=3$ & -1.1656 & 3.0429 & 1.2978 &       & -1.1184 & 3.0196 & 1.3300 &       & -0.6617 & 1.9368 & 1.2243 \\
\cmidrule{1-5}\cmidrule{7-9}\cmidrule{11-13}    \multicolumn{1}{l}{A3} 
& $m_1=1.6$ & -0.1986 & 0.8027 & 0.8737 &       & -0.2654 & 0.8452 & 0.8802 &       & -0.1216 & 0.4168 & 0.6340 \\
& $m_2=3.2$ & -0.2918 & 0.9564 & 0.9334 &       & -0.2704 & 0.9551 & 0.9391 &       & -0.1435 & 0.5343 & 0.7167 \\
& $\mu_1=2$ & 0.0824 & 0.2170 & 0.4585 &       & 0.1353 & 0.2321 & 0.4623 &       & 0.0492 & 0.1331 & 0.3615 \\
& $\mu_2=4$ & 0.2642 & 0.9766 & 0.9523 &       & 0.2180 & 0.9729 & 0.9620 &       & 0.1326 & 0.5618 & 0.7377 \\
& $\beta=1$ & -0.0273 & 0.2381 & 0.4872 &       & -0.0586 & 0.1997 & 0.4430 &       & -0.0337 & 0.0877 & 0.2943 \\
\cmidrule{1-5}\cmidrule{7-9}\cmidrule{11-13}    \multicolumn{1}{l}{A4} 
& $m_1=1.6$ & -0.1950 & 0.8568 & 0.9048 &       & -0.2238 & 0.8574 & 0.8985 &       & -0.1211 & 0.3614 & 0.5889 \\
& $m_2=3.2$ & -0.3265 & 0.9160 & 0.8997 &       & -0.3019 & 0.8859 & 0.8915 &       & -0.1173 & 0.3346 & 0.5664 \\
& $\mu_1=3$ & 0.0973 & 0.5389 & 0.7276 &       & 0.1255 & 0.5382 & 0.7228 &       & 0.0545 & 0.3437 & 0.5837 \\
& $\mu_2=6$ & 0.3677 & 2.3501 & 1.4883 &       & 0.2772 & 2.2471 & 1.4732 &       & 0.1136 & 1.2666 & 1.1197 \\
& $\beta=2$ & -0.1717 & 0.6395 & 0.7810 &       & -0.1821 & 0.6048 & 0.7561 &       & -0.0801 & 0.2873 & 0.5300 \\
\bottomrule
\end{tabular}}
\end{table}

%tabel2
\begin{table}[!t]
    \centering
    \caption{Estimated parameters with $n$=100 for scenarios B1-B4.}\label{table2}
    \resizebox{\linewidth}{!}{
    \begin{tabular}{rlrrrcrrrcrrr}
    \toprule
    &       &  &  \multicolumn{1}{l}{YW}     &       &       &  &  \multicolumn{1}{l}{CLS}     &       &       &  &   \multicolumn{1}{l}{CML}    &  \\
    \cmidrule{3-5}\cmidrule{7-9}\cmidrule{11-13}    \multicolumn{1}{l}{Model} & Para. & \multicolumn{1}{l}{bias} & \multicolumn{1}{l}{MSE} & \multicolumn{1}{l}{SD} &       & \multicolumn{1}{l}{bias} & \multicolumn{1}{l}{MSE} & \multicolumn{1}{l}{SD} &       & \multicolumn{1}{l}{bias} & \multicolumn{1}{l}{MSE} & \multicolumn{1}{l}{SD} \\
    \midrule
    \multicolumn{1}{l}{B1} & $m_1=3.2$ & -0.3600 & 1.0631 & 0.9662 &       & -0.4163 & 1.1289 & 0.9775 &       & -0.1682 & 0.7235 & 0.8338 \\
              & $m_2=4.8$ & -0.6875 & 1.3376 & 0.9300 &       & -0.5740 & 1.1493 & 0.9054 &       & -0.2495 & 0.5882 & 0.7252 \\
              & $\mu_1=2$ & 0.2019 & 0.3426 & 0.5494 &       & 0.2460 & 0.3669 & 0.5535 &       & 0.0861 & 0.1898 & 0.4271 \\
              & $\mu_2=4$ & 0.8867 & 2.3266 & 1.2411 &       & 0.6575 & 1.9404 & 1.2281 &       & 0.2633 & 0.7230 & 0.8085 \\
              & $\phi=1$ & 0.0013 & 5.3856 & 2.3207 &       & -0.0650 & 3.0786 & 1.7534 &       & 0.2418 & 2.0038 & 1.3947 \\
    \cmidrule{1-5}\cmidrule{7-9}\cmidrule{11-13}    \multicolumn{1}{l}{B2} & $m_1=3.2$ & -0.3283 & 1.0105 & 0.9501 &       & -0.3237 & 1.0611 & 0.9779 &       & -0.1411 & 0.6730 & 0.8081 \\
              & $m_2=4.8$ & -0.7225 & 1.2518 & 0.8543 &       & -0.5162 & 0.9634 & 0.8348 &       & -0.2677 & 0.5686 & 0.7049 \\
              & $\mu_1=3$ & 0.2511 & 0.5991 & 0.7322 &       & 0.2382 & 0.6245 & 0.7535 &       & 0.0907 & 0.3300 & 0.5672 \\
              & $\mu_2=6$ & 1.2880 & 4.3934 & 1.6536 &       & 0.7566 & 3.3128 & 1.6554 &       & 0.3738 & 1.4328 & 1.1371 \\
              & $\phi=2$ & -0.0055 & 11.5228 & 3.3945 &       & -0.1837 & 7.3940 & 2.7130 &       & 0.4173 & 4.9265 & 2.1800 \\
    \cmidrule{1-5}\cmidrule{7-9}\cmidrule{11-13}    \multicolumn{1}{l}{B3} & $m_1=1.6$ & -0.2265 & 0.7744 & 0.8503 &       & -0.3460 & 0.8842 & 0.8744 &       & -0.2064 & 0.6958 & 0.8082 \\
              & $m_2=3.2$ & -0.3887 & 1.0603 & 0.9536 &       & -0.3199 & 1.0070 & 0.9511 &       & -0.1692 & 0.5617 & 0.7301 \\
              & $\mu_1=2$ & 0.0844 & 0.1655 & 0.3980 &       & 0.1622 & 0.1951 & 0.4108 &       & 0.0744 & 0.1414 & 0.3686 \\
              & $\mu_2=4$ & 0.3686 & 1.0304 & 0.9458 &       & 0.2752 & 0.9792 & 0.9505 &       & 0.1347 & 0.4324 & 0.6436 \\
              & $\phi=1$& -0.0501 & 0.7463 & 0.8624 &       & -0.0569 & 0.6211 & 0.7860 &       & 0.0241 & 0.4115 & 0.6410 \\
\cmidrule{1-5}\cmidrule{7-9}\cmidrule{11-13}    \multicolumn{1}{l}{B4} 
& $m_1=1.6$ & -0.1740 & 0.7961 & 0.8751 &       & -0.2498 & 0.8495 & 0.8872 &       & -0.1530 & 0.6659 & 0.8015 \\
& $m_2=3.2$ & -0.4189 & 1.0177 & 0.9177 &       & -0.2271 & 0.9283 & 0.9363 &       & -0.1706 & 0.5663 & 0.7329 \\
& $\mu_1=3$ & 0.1049 & 0.3085 & 0.5455 &       & 0.1625 & 0.3305 & 0.5515 &       & 0.0771 & 0.2680 & 0.5119 \\
& $\mu_2=6$ & 0.5561 & 1.8922 & 1.2582 &       & 0.2232 & 1.7260 & 1.2947 &       & 0.1547 & 0.8900 & 0.9306 \\
& $\phi=2$ & -0.0170 & 1.8464 & 1.3587 &       & -0.0481 & 1.5444 & 1.2418 &       & 0.1378 & 1.1480 & 1.0625 \\
\bottomrule
\end{tabular}}
\end{table}

In the simulation studies, the parameters are selected as follows with $l=8$. For I-BRCMNBINAR(1) model, we consider the parameter estimation for the unknown parameter $(m_1,m_2,\mu_1,\mu_2,\beta)$ with the following scenarios:

Scenario A1. $(m_1,m_2,\mu_1,\mu_2,\beta)=(3.2,4.8,2,4,2)$,

Scenario A2. $(m_1,m_2,\mu_1,\mu_2,\beta)=(3.2,4.8,3,6,3)$,

Scenario A3. $(m_1,m_2,\mu_1,\mu_2,\beta)=(1.6,3.2,2,4,1)$,

Scenario A4. $(m_1,m_2,\mu_1,\mu_2,\beta)=(1.6,3.2,3,6,2)$.

%table3
\begin{table}[!t]
\centering
\caption{Estimated parameters $\pmb{\sigma}^2=(\sigma^2_{\alpha_1},\sigma^2_{\alpha_2})^T$ for each scenario.}
\label{table3}
\resizebox{\linewidth}{!}{
\begin{tabular}{rrlrrrrrrr}
\toprule
&       &       &       & \multicolumn{1}{l}{CLS} &       &       &       & \multicolumn{1}{l}{CML} &  \\
\cmidrule{5-7}\cmidrule{9-10} \multicolumn{1}{l}{Series} & \multicolumn{1}{l}{$n$} & Para. &       & \multicolumn{1}{l}{bias} & \multicolumn{1}{l}{SD} & \multicolumn{1}{l}{Per.} &      & \multicolumn{1}{l}{bias} & \multicolumn{1}{l}{SD} \\
\midrule
\multicolumn{1}{l}{A1} & 100   & $\sigma^2_{\alpha_1}=0.24$ &       & -0.0566  & 0.0451  & 0.6700  &       & -0.0055  & 0.0191  \\
            &       & $\sigma^2_{\alpha_2}=0.4971$ &       & -0.0851  & 0.0504  & 0.8680  &       & -0.0055  & 0.0233  \\
            & 200   & $\sigma^2_{\alpha_1}=0.24$ &       & -0.0429  & 0.0404  & 0.8300  &       & -0.0029  & 0.0140  \\
            &       & $\sigma^2_{\alpha_2}=0.4971$ &       & -0.0695  & 0.0434  & 0.9750  &       & -0.0038  & 0.0162  \\
            & 300   & $\sigma^2_{\alpha_1}=0.24$ &       & -0.0366  & 0.0374  & 0.9040  &       & -0.0021  & 0.0116  \\
            &       & $\sigma^2_{\alpha_2}=0.4971$ &       & -0.0637  & 0.0390  & 0.9960  &       & -0.0030  & 0.0136  \\
            & 400   & $\sigma^2_{\alpha_1}=0.24$ &       & -0.0332  & 0.0352  & 0.9380  &       & -0.0017  & 0.0104  \\
            &       & $\sigma^2_{\alpha_2}=0.4971$ &       & -0.0603  & 0.0388  & 0.9990  &       & -0.0024  & 0.0116  \\
\midrule
\multicolumn{1}{l}{A2} & 100   & $\sigma^2_{\alpha_1}=0.24$ &       & -0.0592  & 0.0484  & 0.6150  &       & -0.0160  & 0.0309  \\
            &       & $\sigma^2_{\alpha_2}=0.4971$ &       & -0.0860  & 0.0522  & 0.8510  &       & -0.0251  & 0.0495  \\
            & 200   & $\sigma^2_{\alpha_1}=0.24$ &       & -0.0429  & 0.0434  & 0.7980  &       & -0.0039  & 0.0126  \\
            &       & $\sigma^2_{\alpha_2}=0.4971$ &       & -0.0705  & 0.0453  & 0.9580  &       & -0.0050  & 0.0145  \\
            & 300   & $\sigma^2_{\alpha_1}=0.24$ &       & -0.0365  & 0.0389  & 0.8880  &       & -0.0034  & 0.0105  \\
            &       & $\sigma^2_{\alpha_2}=0.4971$ &       & -0.0630  & 0.0413  & 0.9830  &       & -0.0049  & 0.0135  \\
            & 400   & $\sigma^2_{\alpha_1}=0.24$ &       & -0.0329  & 0.0374  & 0.9110  &       & -0.0031  & 0.0091  \\
            &       & $\sigma^2_{\alpha_2}=0.4971$ &       & -0.0581  & 0.0396  & 0.9990  &       & -0.0044  & 0.0113  \\
\midrule
\multicolumn{1}{l}{A3} & 100   & $\sigma^2_{\alpha_1}=0.0743$
 &       & -0.0254  & 0.0379  & 0.5440  &       & -0.0021  & 0.0155  \\
            &       & $\sigma^2_{\alpha_2}=0.24$ &       & -0.0500  & 0.0431  & 0.7360  &       & -0.0034  & 0.0224  \\
            & 200   & $\sigma^2_{\alpha_1}=0.0743$ &       & -0.0169  & 0.0343  & 0.6530  &       & -0.0011  & 0.0109  \\
            &       & $\sigma^2_{\alpha_2}=0.24$ &       & -0.0378  & 0.0372  & 0.8830  &       & -0.0019  & 0.0153  \\
            & 300   & $\sigma^2_{\alpha_1}=0.0743$ &       & -0.0135  & 0.0310  & 0.7460  &       & -0.0007  & 0.0093  \\
            &       & $\sigma^2_{\alpha_2}=0.24$ &       & -0.0315  & 0.0342  & 0.9440  &       & -0.0011  & 0.0120  \\
            & 400   & $\sigma^2_{\alpha_1}=0.0743$ &       & -0.0114  & 0.0286  & 0.7970  &       & -0.0004  & 0.0081  \\
            &       & $\sigma^2_{\alpha_2}=0.24$ &       & -0.0285  & 0.0317  & 0.9720  &       & -0.0009  & 0.0106  \\
\midrule
\multicolumn{1}{l}{A4} & 100   & $\sigma^2_{\alpha_1}=0.0743$ &       & -0.0282  & 0.0452  & 0.4620  &       & -0.0022  & 0.0138  \\
            &       & $\sigma^2_{\alpha_2}=0.24$ &       & -0.0555  & 0.0475  & 0.6550  &       & -0.0030  & 0.0178  \\
            & 200   & $\sigma^2_{\alpha_1}=0.0743$ &       & -0.0187  & 0.0418  & 0.5770  &       & -0.0013  & 0.0090  \\
            &       & $\sigma^2_{\alpha_2}=0.24$ &       & -0.0413  & 0.0411  & 0.8350  &       & -0.0011  & 0.0125  \\
            & 300   & $\sigma^2_{\alpha_1}=0.0743$ &       & -0.0137  & 0.0385  & 0.6570  &       & -0.0009  & 0.0068  \\
            &       & $\sigma^2_{\alpha_2}=0.24$ &       & -0.0347  & 0.0377  & 0.9030  &       & -0.0010  & 0.0097  \\
            & 400   & $\sigma^2_{\alpha_1}=0.0743$ &       & -0.0115  & 0.0356  & 0.7130  &       & -0.0008  & 0.0060  \\
            &       & $\sigma^2_{\alpha_2}=0.24$ &       & -0.0303  & 0.0357  & 0.9380  &       & -0.0009  & 0.0084  \\
      \bottomrule
      \end{tabular}}%
\end{table}

For II-BRCMNBINAR(1) model, we focus on the estimation of the unknown parameter $(m_1,m_2,\mu_1,\mu_2,\phi)$ with the following scenarios:

Scenario B1. $(m_1,m_2,\mu_1,\mu_2,\phi)=(3.2,4.8,2,4,1)$,

Scenario B2. $(m_1,m_2,\mu_1,\mu_2,\phi)=(3.2,4.8,3,6,2)$,

Scenario B3. $(m_1,m_2,\mu_1,\mu_2,\phi)=(1.6,3.2,2,4,1)$, %B3与B4换序了

Scenario B4. $(m_1,m_2,\mu_1,\mu_2,\phi)=(1.6,3.2,3,6,2)$.

We consider the parameter estimation problems using YW, CLS and CML methods respectively. Let $\textbf{\textit{X}}_1,\textbf{\textit{X}}_2,\cdots,\textbf{\textit{X}}_n$ be the data generated by I-BRCMNBINAR(1) and II-BRCMNBINAR(1) models with corresponding true parameters. The simulation results are summarized in Tables \ref{table1}-\ref{table2}. All the calculations are performed under $\mathcal{R}$ software with 1000 replications.

We compute the empirical bias, mean squared error (MSE) and standard deviation (SD) for each combination of the parameters to evaluate the efficiency of these three methods. From Tables \ref{table1}-\ref{table2}, we find that the results of YW-estimator and CLS-estimator are closed. %, which is consistent with the fact that they have the same asymptotic distribution
CLS-estimator is slightly better than YW-estimator in terms of MSE. From the simulation results, the CML-estimator has smaller bias, MSE and SD, and better robustness than YW-estimator and CLS-estimator. For the bivariate time series models, the time required for CML estimation is significantly longer than the other two methods, especially as the sample size increases. For all series in Tables \ref{table1}-\ref{table2}, the bias, MSE and SD of these three estimators decrease as the sample size increases for all cases. It implies that the estimators are consistent for all parameters, and these three estimation methods are reliable.

Figs.\ref{figure1}-\ref{figure4} display the boxplots of the estimates for scenarios A1, A4, B1 and B4 respectively. As we can see from the figures, median of these estimators are closer to the real parameter values as the sample size increases. Thus, the figure illustrates that the performance of all these three estimation methods increase gradually. Moreover, by comparing the results in the graphs, we can see that the interquartile ranges and the overall range of the estimators both become narrower, which indicates low dispersion. In addition, the performance of CML-estimator is significantly better than other two methods in terms of the location with median estimates that are apparently closer to the real values of the parameter.

We consider the estimates of the variance $\pmb{\sigma}^2=(\sigma^2_{\alpha_1},\sigma^2_{\alpha_2})^T$, where the innovations follow the BVNB distribution as an example. The results of the simulation study are listed in Table \ref{table3}. Notice that there is a positive probability that the CLS-estimator of $\pmb{\sigma}^2$ is negative with finite sample size and hence we provide the percentage of positive estimates for $\pmb{\sigma}^2$ based on two-step CLS method in the Per. column of Table \ref{table3}. From Table \ref{table3}, we find that the percentage of positive estimates for $\pmb{\sigma}^2$ increases with respect to the sample size $n$. More specifically, the percentage of negative estimates for $\pmb{\sigma}^2$ is smaller for larger variances $\pmb{\sigma}^2$. In general, two-step CLS and CML methods both can bring good estimators of the parameter $\pmb{\sigma}^2$ for larger sample size.
\begin{figure}[!t]
\centering
\includegraphics[width=14cm,height=10cm]{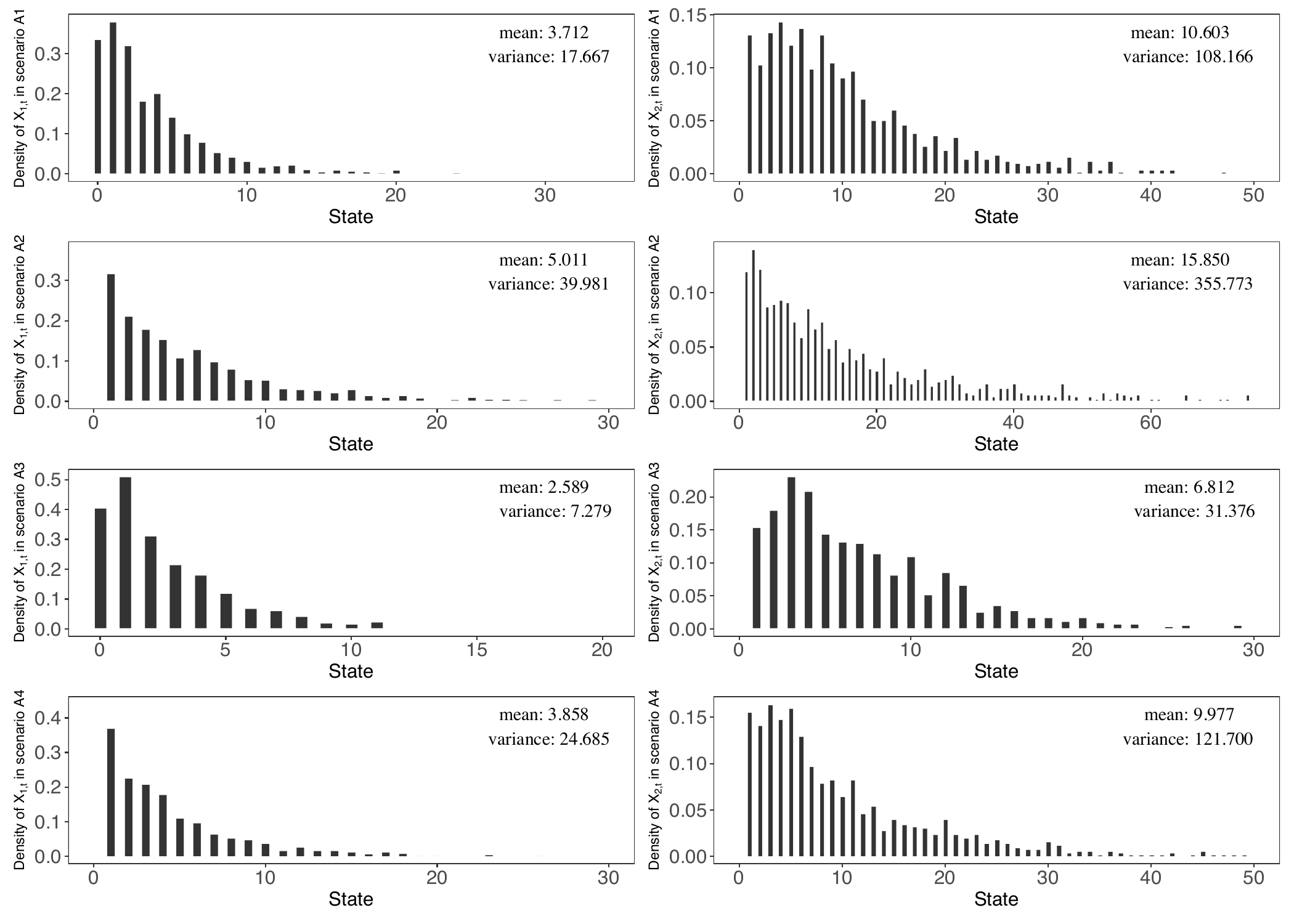}
\caption{Stationary marginal distributions for Scenarios A1-A4.}
\label{figure7}
\end{figure}

\section{Coherent forecasting for BRCMNBINAR(1) model}\label{SEC5}
One of the main applications of the BINAR(1) process is to predict $\textbf{\textit{X}}_{n+h}$ for $h\geqslant 1$ based on past observations $\textbf{\textit{X}}_1,\textbf{\textit{X}}_2,\cdots,\textbf{\textit{X}}_n$. One common approach for time series forecasting is to use the conditional expectation, which yields forecasts with minimizing the mean square error. However, this method is unsatisfactory sometimes as it might lead to a non-integer value. To overcome the disadvantage of point forecast, some scholars suggested employing the coherent forecasting techniques, which will only produce forecasts on $\mathbb{N}_0$. For example, Freeland and McCabe (2004) considered the median of the $h$-step-ahead conditional forecast distribution. Furthermore, Li et al. (2023) and Yang et al. (2023b) generalized the coherent forecasting method to the covariate-driven INAR models.

%Bu and McCabe (2008) and Wei\ss (2010), 
To forecast the BRCMNBINAR(1) process, we use the forecasting distributions of $\textbf{\textit{X}}_{n+h}$ given on the observation $\textbf{\textit{X}}_{n}$ over all horizons $h\in\mathbb{N}$. For a Markov chain with finite states, the $h$-step-ahead conditional distribution of $\textit{\textbf{X}}_{n+h}$ given on $\textit{\textbf{X}}_{n}$ is expressed as $P(\textit{\textbf{X}}_{n+h}=\textit{\textbf{x}}_{n+h}|\textit{\textbf{X}}_{n}=\textit{\textbf{x}}_{n})=[\textbf{\textit{P}}^h]_{\textit{\textbf{x}}_{n+h},\textit{\textbf{x}}_{n}}$, where $\textit{\textbf{P}}$ is the transition matrix with the elements $p_{\textit{\textbf{k}}|\textit{\textbf{s}}}=P(\textit{\textbf{X}}_{t}=\textit{\textbf{k}}|\textit{\textbf{X}}_{t-1}=\textit{\textbf{s}})$ defined in (\ref{eq3}). In principle, the BRCMNBINAR(1) process should take infinite values on state space $\mathbb{N}_0^2$, which makes it difficult to compute the marginal distribution. In fact, as discussed by Yang et al. (2023a), we can choose two sufficiently large positive constants $M_1$ and $M_2$, such that the probability of $\textit{\textbf{X}}_{t}$ larger than $(M_1, M_2)$ is negligible. The next question is how to choose $M_1$ and $M_2$ in practice. As suggested by Yang et al. (2023b), we can choose an integer $M=\max\{M_1,M_2\}$ for simplicity. Therefore, the state space is given by 
\begin{align*}
\mathbb{S}=\{(0,0),(0,1),\cdots,(0,M),\;(1,0)\cdots,(1,M),\; (M,0)\cdots,(M,M)\}.
\end{align*}
Let ${\pmb{\mit{P}}}_M$ denote the $(M+1)^2\times(M+1)^2$ transition matrix of the BRCMNBINAR(1) model which admits the form
\begin{align*}
\textbf{\textit{P}}_{M}=
\begin{pmatrix}
p_{0,0|0,0} & p_{0,0|0,1} &\cdots & p_{0,0|M,M}\\
p_{0,1|0,0} & p_{0,1|0,1} &\cdots & p_{0,1|M,M}\\
\cdots&\cdots&\cdots&\cdots\\
p_{M,M|0,0} & p_{M,M|0,1} &\cdots & p_{M,M|M,M}
\end{pmatrix}.
\end{align*}
The approximated marginal distribution $\pmb{\pi}_M=(\pi_{0,0}, \pi_{0,1}, \pi_{0,M},...,\pi_{1,M}, \pi_{2,M}, \pi_{M,M})^T$ can be obtained by solving the equation $\pmb{\mit P}_M\pmb{\pi}_M=\pmb{\pi}_M$. Then the $h$-step-ahead forecasting conditional distribution based on maximum likelihood estimator $p_h(\textbf{\textit{X}}_{n+h}|\textbf{\textit{X}}_{n},\hat{\pmb\tau}^{CML})$ can be obtained approximatively as 
\begin{align*}
p_h(\textit{\textbf{X}}_{n+h}=\textit{\textbf{x}}_{n+h}|\textit{\textbf{X}}_{n}=\textit{\textbf{x}}_{n},\hat{\pmb\tau}^{CML})=P(\textbf{\textit{X}}_{n+h}=\textbf{\textit{x}}_{n+h}|\textbf{\textit{X}}_{n}=\textbf{\textit{x}}_{n})=[\textbf{\textit{P}}_{M}^h]_{\textit{\textbf{x}}_{n+h},\textit{\textbf{x}}_{n}}.
\end{align*}

Under standard regularity conditions, the maximum likelihood estimator $\hat{\pmb\tau}^{CML}$ is asymptotically normal distributed around the true value $\pmb\tau$, i.e. $\sqrt{n-1}(\hat{\pmb\tau}^{CML}-\pmb\tau)\xrightarrow{L} N(\textbf{\textit{0}},\textbf{\textit{I}}^{-1})$, where $\textbf{\textit{I}}$ is the Fisher information matrix. In the following theorem, we work out the asymptotic distribution and the confidence interval for the one-step-ahead prediction of $p_h(\textbf{\textit{X}}_{n+h}|\textbf{\textit{X}}_{n},\hat{\pmb\tau}^{CML})$.
\begin{theorem}
For fixed $\textit{\textbf{X}}_1\in\mathbb{N}_0^2$, the quantity $p_h(\textit{\textbf{X}}_{n+h}|\textit{\textbf{X}}_n,\hat{\pmb\tau}^{CML})$ has an asymptotically normal distribution, that is 
\begin{align*}
\sqrt{n-1}\big(p_h(\textit{\textbf{X}}_{n+h}|\textit{\textbf{X}}_{n},\hat{\pmb\tau}^{CML})-p_h(\textit{\textbf{X}}_{n+h}|\textit{\textbf{X}}_{n},\pmb\tau)\big)\xrightarrow{L}N(0,\sigma^2),
\end{align*} 
where $\sigma^2=\textit{\textbf{D}}\textit{\textbf{I}}^{-1}\textit{\textbf{D}}^T$, $\textit{\textbf{D}}=\partial p_h(\textit{\textbf{X}}_{n+h}|\textit{\textbf{X}}_{n},\pmb\tau)/\partial{\pmb\tau}^T$ is a vector of partial derivatives. Furthermore, the $100(1-\alpha)\%$ confidence interval for $p_h(\textit{\textbf{X}}_{n+h}|\textit{\textbf{X}}_{n},\hat{\pmb\tau}^{CML})$ is given by
\begin{align*}
(p_h(\textit{\textbf{X}}_{n+h}|\textit{\textbf{X}}_{n},\hat{\pmb\tau}^{CML})-\frac{\sigma}{\sqrt{n-1}}u_{1-\frac{\alpha}{2}}, p_h(\textit{\textbf{X}}_{n+h}|\textit{\textbf{X}}_{n},\ \hat{\pmb\tau}^{CML})+\frac{\sigma}{\sqrt{n-1}}u_{1-\frac{\alpha}{2}}),
\end{align*} 
where $u_{1-\frac{\alpha}{2}}$ is the $(1-\frac{\alpha}{2})$ upper quantile of standard normal distribution.
\end{theorem}

\begin{figure}[!t]
\centering
\includegraphics[width=14cm,height=8cm]{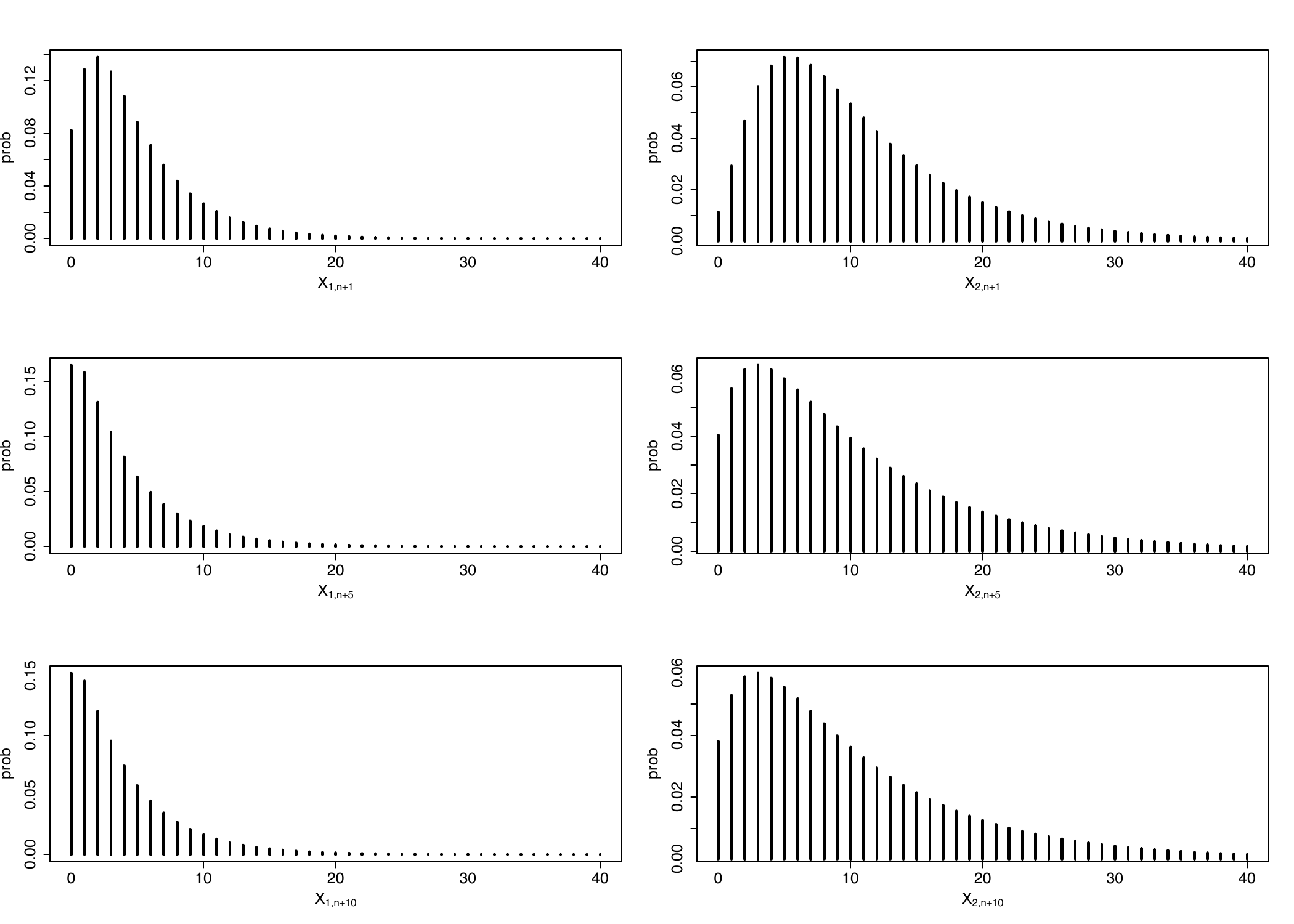}
\caption{Marginal forecast distributions for Scenario A1 with $M=40$ conditional on $\textit{\textbf{X}}_n=(6,10)^T$.}
\label{forecast1}
\end{figure}

\begin{figure}[!t]
\includegraphics[width=14cm,height=8cm]{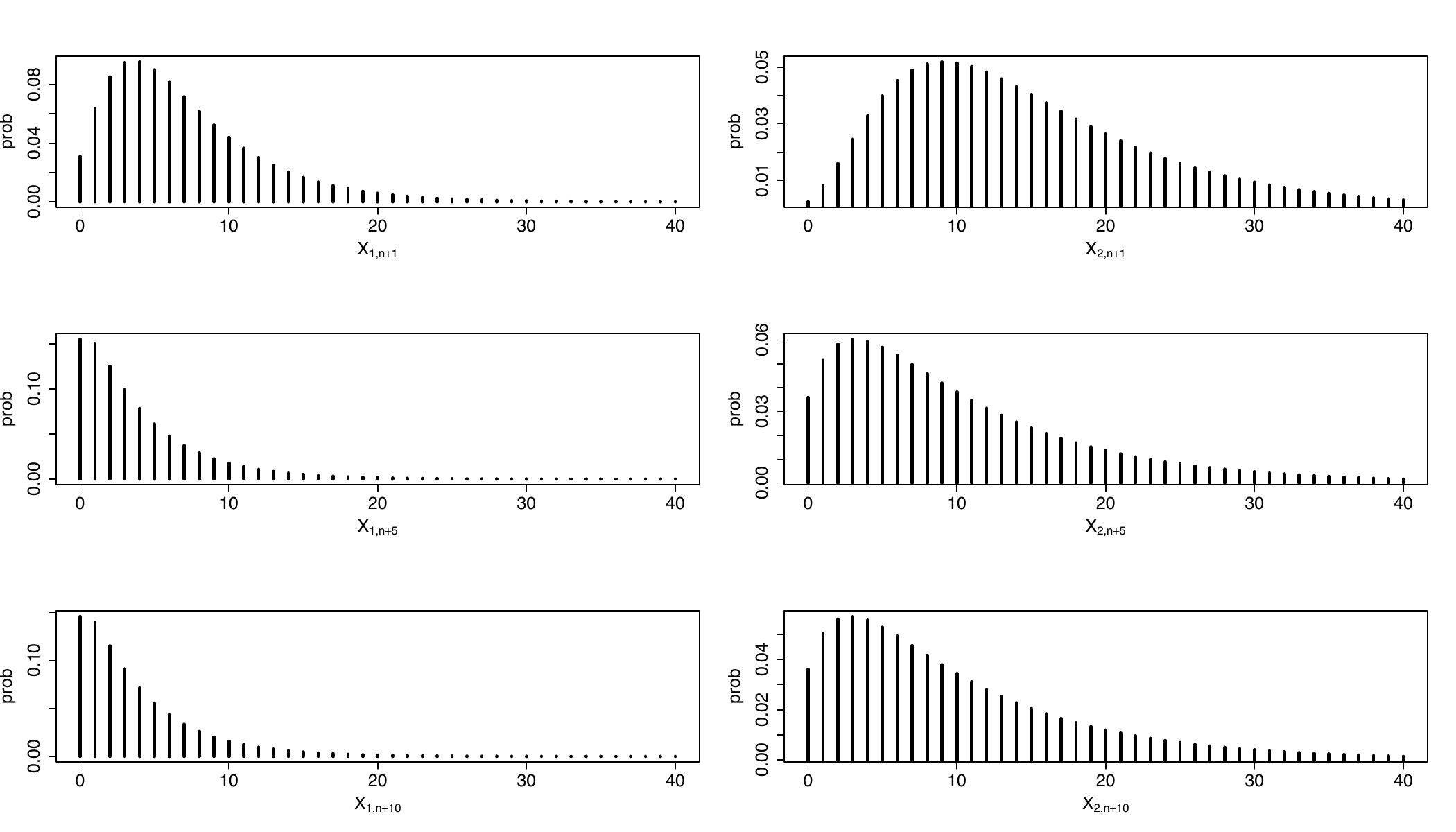}
\caption{Marginal forecast distributions for Scenario A1 with $M=40$ conditional on $\textit{\textbf{X}}_n=(12,18)^T$.}
\label{forecast2}
\end{figure}

From Fig. \ref{figure7}, we can see that the variance is greater than the mean for all scenarios, which indicates that they are overdispersed. The probabilities tend to zero when the states become larger, meaning that the approximated method described in Section \ref{SEC5} is feasible. As an illustration, Figs. \ref{forecast1} and \ref{forecast2} provide the forecasting distributions for the horizons $h=1,5,10$ steps ahead for Scenario A1, i.e. $(m_1,m_2,\mu_1,\mu_2,\beta)=(3.2,4.8,2,4,2)$, conditional on $\textit{\textbf{X}}_n=(6,10)^T$ and $\textit{\textbf{X}}_n=(12,18)^T$ respectively. These conditional marginal distributions converge to the stationary marginal distributions shown in Fig. \ref{figure7} for growing $h$, which is expected from the ergodicity of the BRCMNBINAR(1) process. As in Figs. \ref{forecast1} and \ref{forecast2}, we can see that the figures have similar “shape” of the conditional marginal distributions even if we change the initial observation $\textit{\textbf{X}}_n$, which is consistent with the stationarity of the BRCMNBINAR(1) process. 
%检验准确性？

\section{Real data example}\label{SEC6}
%这个不错
%Richmond Valley	Richmond Valley
%Sexual offences	Drug offences
%Sexual touching, sexual act and other sexual offences	Possession and/or use of cannabis
%=======================
In this section, we apply our proposed BRCMNBINAR(1) process to fit a bivariate count time series of the monthly crime datasets from the NSW Bureau of Crime Statistics and Research. The datasets cover from January 1995 to December 2022. Each dataset comprises 336 observations, categorized by type of offence, month and local government area. In this work, we mainly consider two types of crime: domestic-violence-related assault and drug offences in Bellingen. Assault is the sum of three subcategories: domestic-violence-related assault, non-domestic violence related assault and assault police, where the first category accounts for 40\% of the assault offence. Thus, we mainly consider the subcategory of `domestic-violence-related assault' and label it as `assault'. Aside from the domestic-violence-related assault, we note that drug offences that are recorded by the NSW Police Force include five subcategories: possession drugs, dealing drugs, cultivating cannabis, manufacturing drugs and importing drugs. We observe that most of the counts in each subcategory for each month are zeros and over half of the drug offences are dominated by the subcategory of `possession and/or use of cannabis' during this time period. Thus, we handle this dataset and label it as `drug offence'. In this example, we focus on whether drug offence does affect the domestic-violence-related assault.
\begin{figure}[!t]
\centering
\includegraphics[width=14cm,height=8cm]{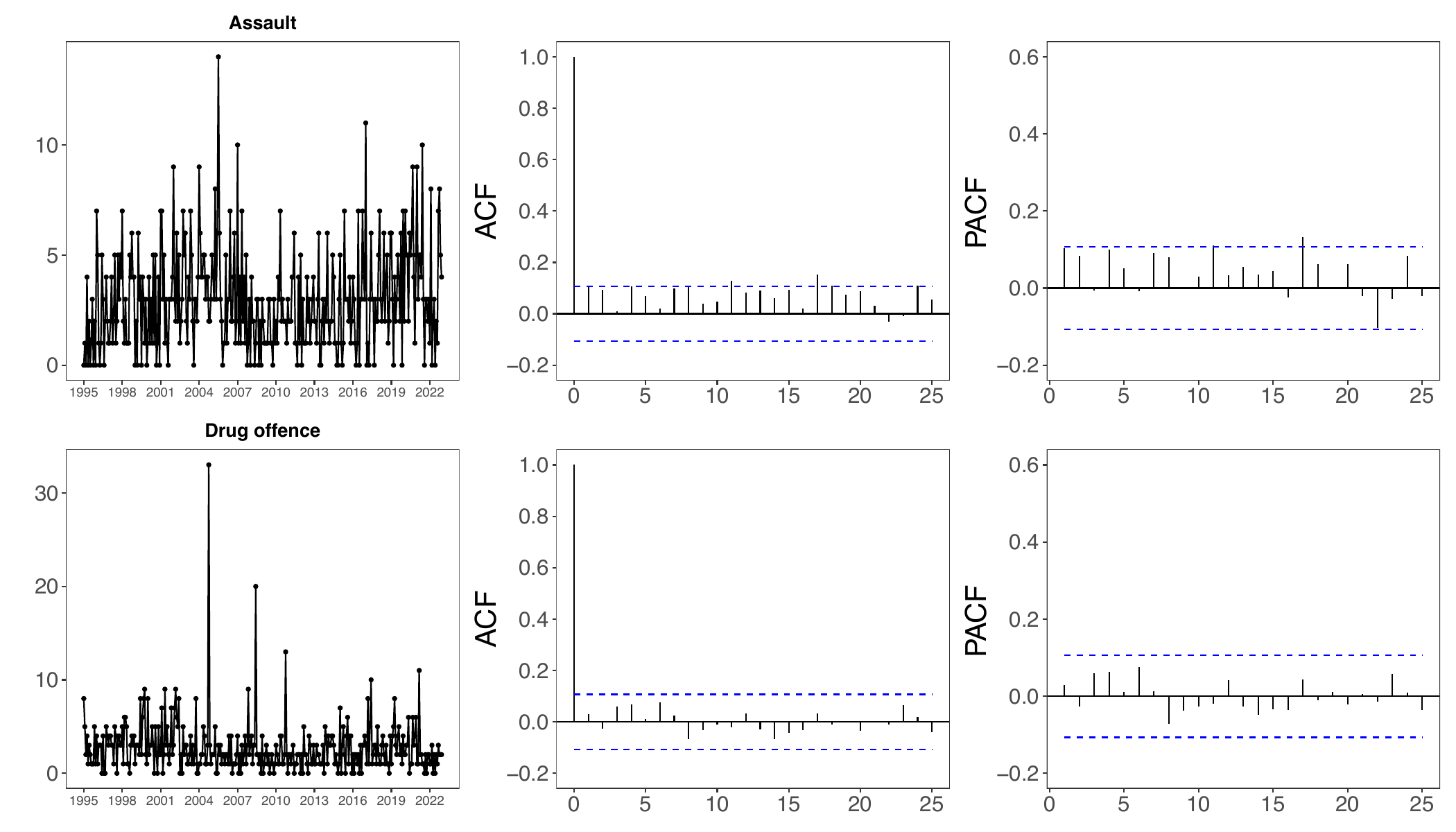}
\includegraphics[width=5cm,height=4.5cm]{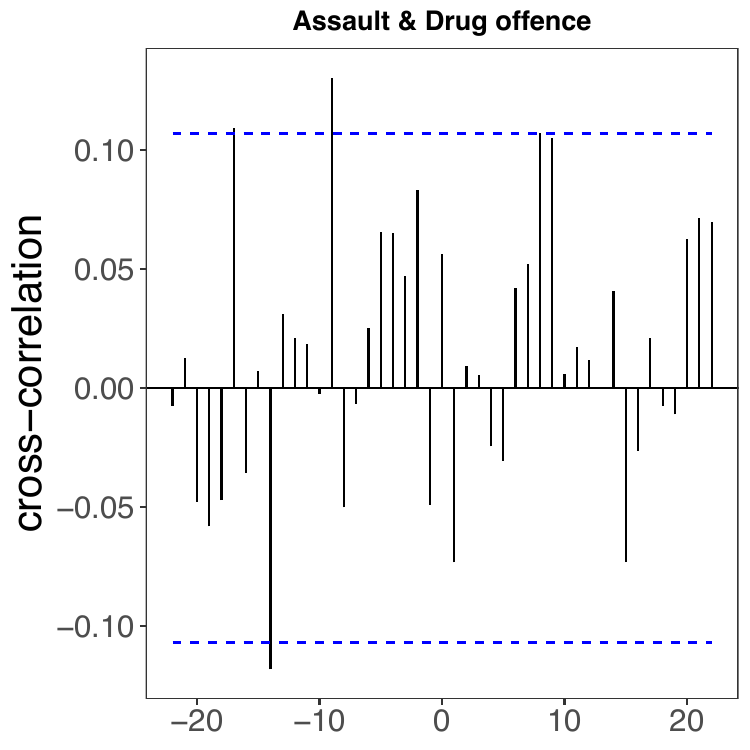}
\caption{Sample path, ACF, PACF and CCF plots of the monthly counts for assault and drug offence in Bellingen.}
\label{dataset_samplepath}
\end{figure}

\begin{figure}[!t]
\centering
\includegraphics[width=14cm,height=5cm]{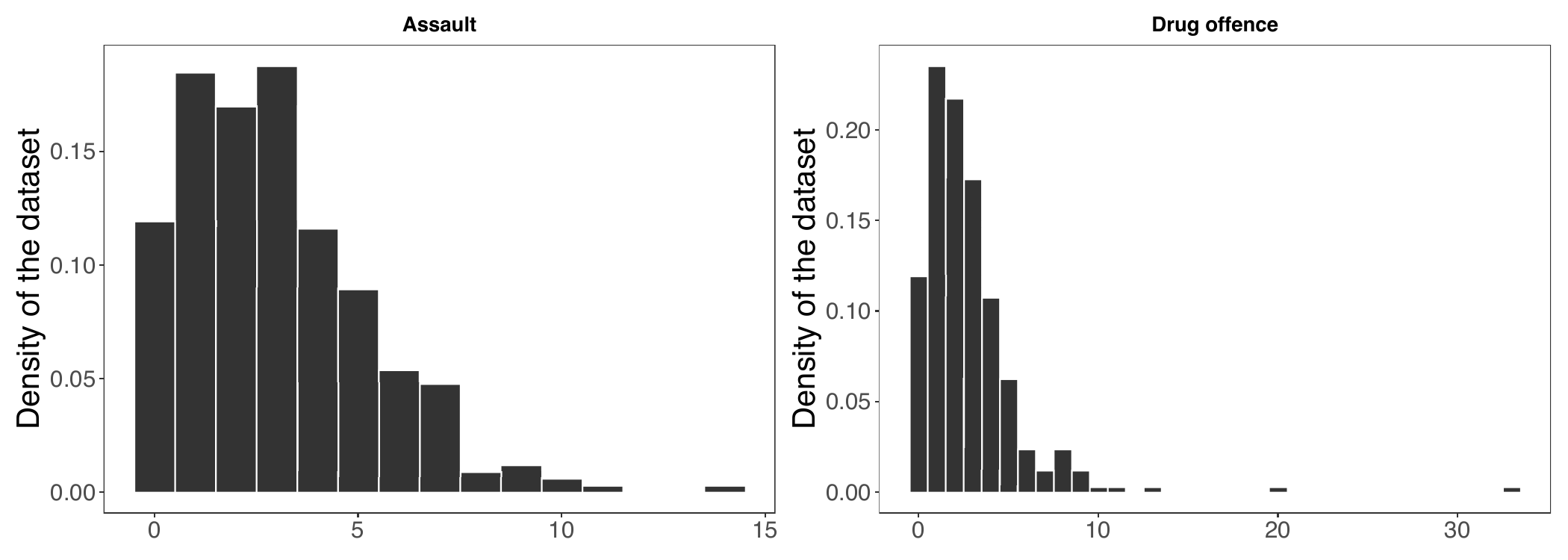}
\caption{Marginal distributions for the series of assault and drug offence.}
\label{dataset_marginal}
\centering
\end{figure}

The sample mean values for assault and drug offence are 2.964 and 2.711, while the sample variances are 5.204 and 8.128 respectively. These bivariate crime series are all overdispersed, and the range of them are [0,14] and [0,33] respectively. Moreover, the sample cross-correlation coefficient between the counts of the bivariate time series is 0.056. The lag 1 autocorrelation coefficient are 0.104 and 0.03 respectively. Therefore, the real data may come from a bivariate integer-valued AR(1) process. Fig. \ref{dataset_samplepath} shows the sample paths, the sample autocorrelation functions (ACF), partial auto-correlation functions (PACF) and the cross-correlation function (CCF) plots of the monthly time series for the assault and drug offence. Regarding the marginal distributions in Fig. \ref{dataset_marginal}, both of them are skewed to the right and exhibit considerable overdispersion characteristics.

We use I-BRCMNBINAR(1) and II-BRCMNBINAR(1) models to fit the datasets and compare them with the following four different models:

$\bullet$ BINAR(1) process with BVNB innovations (Pedeli and Karlis (2011)).

$\bullet$ BINAR(1) process with BP innovations (Pedeli and Karlis (2013a)).

$\bullet$ BNBINAR(1) with BP and BVNB innovations (Zhang et al. (2020)).

%$\bullet$ BRCINAR(1) models with Beta distribution (Yu et al. (2020)).%and Gaussian mixing 

%\textcolor{blue}{$\bullet$ BTINAR(1) process with BP innovations (Yang et al. (2023a))}.

The reason for choosing the above models is that they are based on diagonal-type random matricial operation with binomial thinning operator and modified negative binomial operator. These models are similar to our proposed model. For these above models, we use the CML method to estimate the unknown parameters with $l=5$ and take the Akaike information criterion (AIC) and the Bayesian information criterion (BIC) as the goodness of fit criteria. These two measures illustrate how well the proposed distribution of the model fits the dataset, rather than which model is the best. Table \ref{table4} summarizes the fitting results, including the parameter estimators, standard errors as well as the values for the two goodness-of-fit measures. 
%Motor vehicle theft & drug in bell 还行
%Willoughby Sexual offences & Possession and/or use of cannabis
\begin{table}[t!]
\centering
\caption{CML estimates, SE, AIC and BIC for the offence data.}
\label{table4}
\resizebox{\linewidth}{!}{
\begin{tabular}{cccccc}
\toprule
Model & Para. & CML   & SE    & AIC   & BIC \\
\midrule
BNBINAR(1) with BP innovations & 
$\alpha_1$ & 0.2194 & 0.0534 & 2975.139 & 2994.225 \\
& $\alpha_2$ & 0.3368 & 0.0435 &       &  \\
& $\lambda_1$ & 2.1036  & 0.2188 &       &  \\
& $\lambda_2$ & 1.4447 & 0.1594  &       &  \\
& $\phi$   & 0.3928 & 0.1163 &       &  \\
\midrule
BNBINAR(1) with BVNB innovations 
& $\alpha_1$ & 0.1967 &0.0585 &2870.648 &2889.734\\
& $\alpha_2$ & 0.0236 & 0.0345 &       &  \\
& $\lambda_1$ & 2.1939 & 0.2483 &       &  \\
& $\lambda_2$ & 2.6085 & 0.1727  &       &  \\
& $\beta$  & 0.2912 & 0.0509 &       &  \\
\midrule
I-BRCMNBINAR(1)  %check
& $m_1$ & 0.6784 & 0.2738 & 2800.981&  2820.066 \\
& $m_2$ & 0.9409 & 0.2512 &       &  \\
& $\mu_1$ & 2.4402 & 0.2199 &       &  \\
& $\mu_2$ & 1.8959 & 0.1764 &       &  \\
&  $\beta$ & 0.1949 & 0.0434 &       &  \\
\midrule
II-BRCMNBINAR(1)  
& $m_1$ & 1.1663& 0.2606 & 2822.457 &2841.543 \\
& $m_2$ & 1.1286& 0.2258 &       &  \\
& $\mu_1$   & 2.067& 0.1832&       &  \\
& $\mu_2$   & 1.7472 & 0.1439  &       &  \\
& $\phi$   & 0.4882 & 0.1349 &       &  \\
\midrule
BINAR(1) process with BVNB innovations  & 
$\alpha_1$ & 0.4138 & 0.0237 & 4097.711 & 4116.797\\
& $\alpha_2$ &  0.3546 & 0.0257  &  & \\
& $\lambda_1$ & 1.7420 & 0.1238 &  &\\
& $\lambda_2$ & 1.7383 & 0.1243    &  &\\
& $\beta$ & 0.8089 & 0.1132    &  &\\
\midrule
BINAR(1) process with BP innovations  & 
$\alpha_1$ & 0.0732 & 0.0352 & 3028.814 & 3047.899\\
& $\alpha_2$ &  0.0210 & 0.0260  &  & \\
& $\lambda_1$ & 2.7559 & 0.1359 &  &\\
& $\lambda_2$ & 2.6388 & 0.1127  &  &\\
& $\phi$ & 0.1605 & 0.0975    &  & \\
\bottomrule
\end{tabular}}
\end{table}

\begin{figure}[!t]
\centering
\includegraphics[width=14cm,height=6cm]{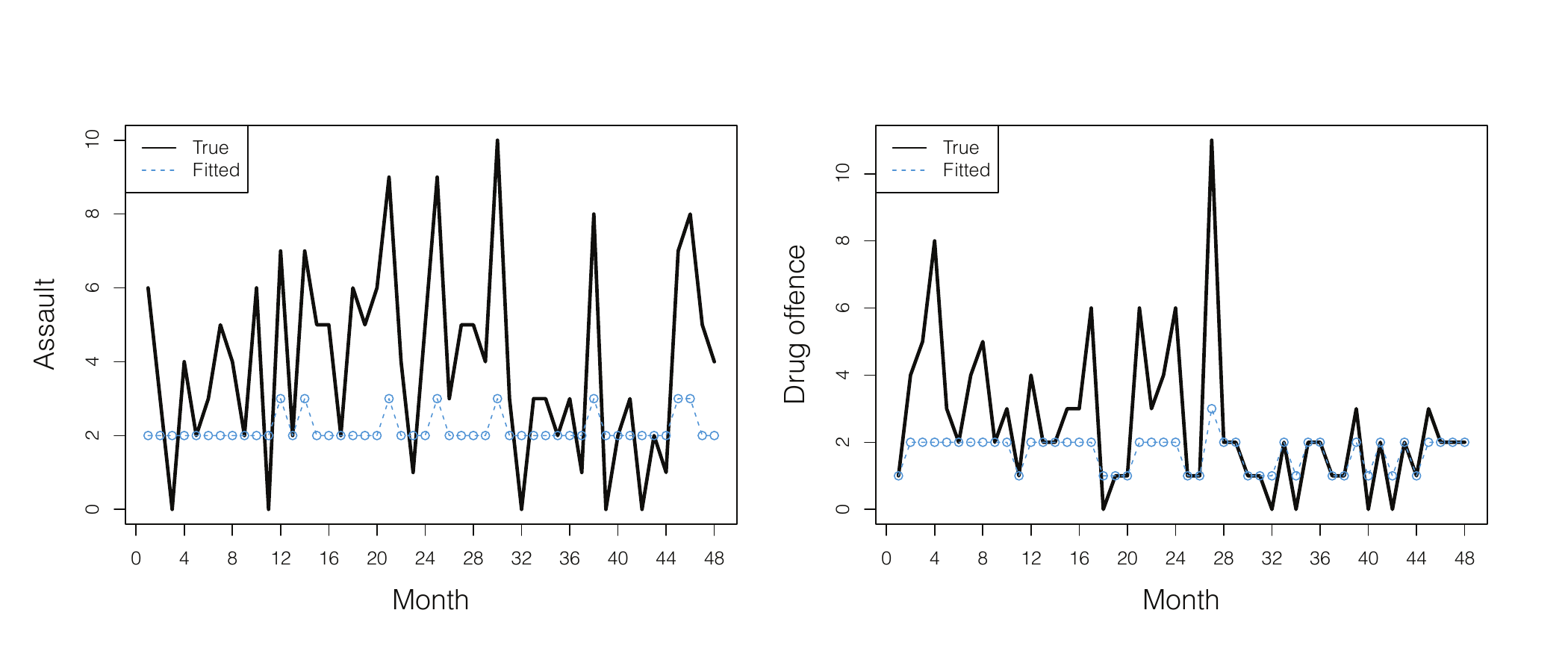}
\caption{Observations of the series of assault and drug offence and the corresponding one-step-ahead forecasts.}
\label{dataset_fitted}
\centering
\includegraphics[width=14cm,height=5cm]{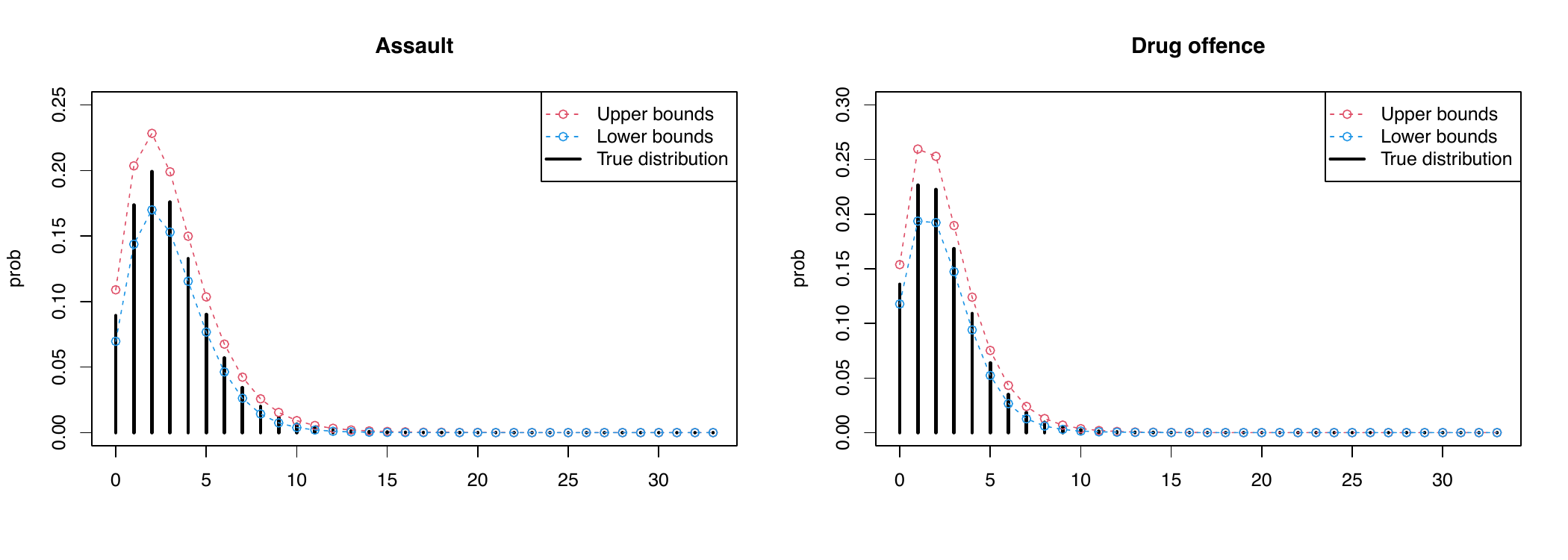}
\caption{One-step-ahead marginal forecast distribution for the series of assault and drug offence conditional on $\textbf{\textit{X}}_{288}=(6,4)^T$ and their 95\% confidence intervals.}
\label{dataset_CI}
\end{figure}

%first 288 
Furthermore, we use the method stated in Section \ref{SEC5} to produce coherent forecasts for the offence data with I-BRCMNBINAR(1) model. We divide the offence data of assault and drug offence into two parts. We use the data from January 1995 to December 2017 to estimate the unknown parameters, and leave the remaining data from January 2018 to December 2022 for the one-step ahead forecast. Fig. \ref{dataset_fitted} provides the corresponding marginal forecasting distributions of the assault and drug offence data given that $\textbf{\textit{X}}_{288}=(6,4)^T$. Fig. \ref{dataset_CI} shows the 95\% confidence intervals and concludes that the most possible one-step-ahead predictive values are equal to 2 and 1 for the assault and drug offence series respectively.

\begin{figure}[!t]
\centering
\includegraphics[width=14cm,height=4cm]{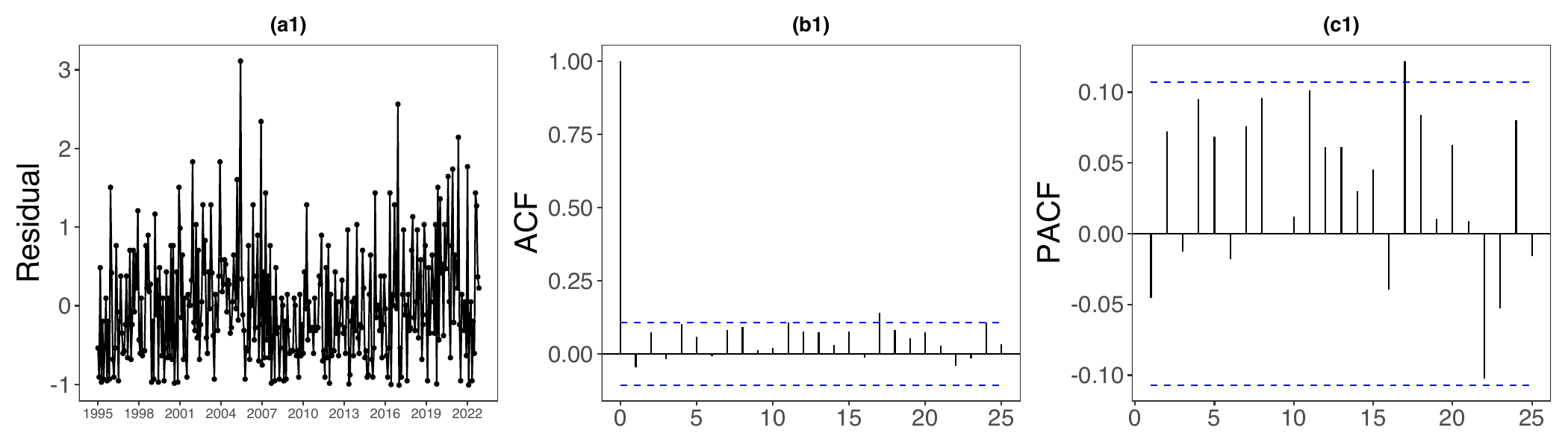}
\includegraphics[width=14cm,height=4cm]{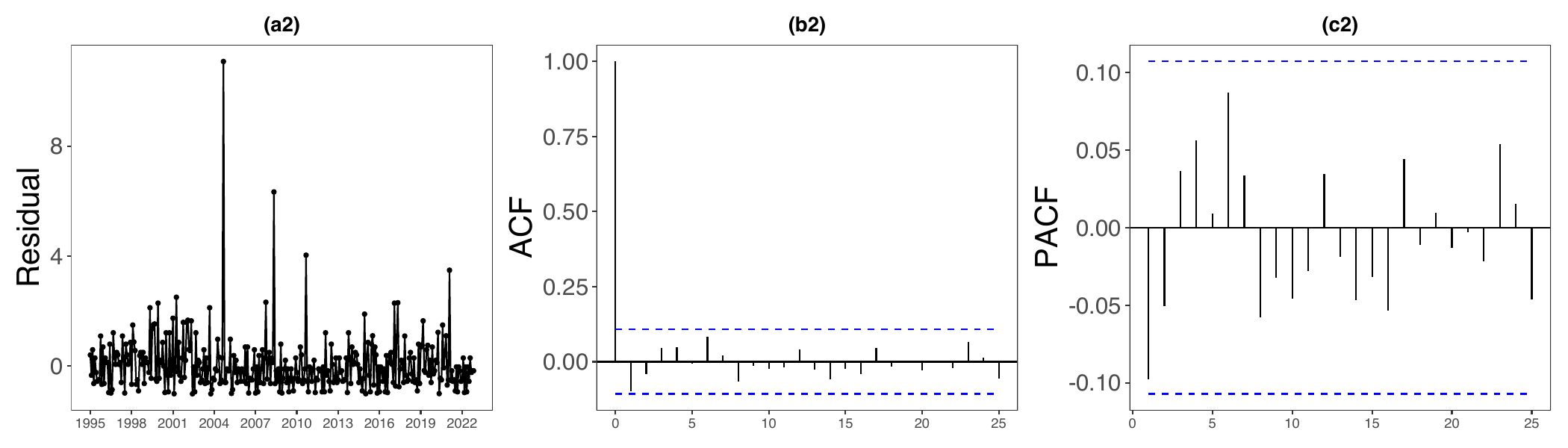}
\caption{Diagnostic checking plots of the fitted I-BRCMNBINAR(1) model. (a1)-(a2) The trace plots of standard Pearson residuals; (b1)-(b2) ACF plots of standard Pearson residuals; (c1)-(c2) PACF plots of standard Pearson residuals.}
\label{Pearson_residual}
\end{figure}

In the following, we conduct the diagnostic checking for the fitted I-BRCMNBINAR(1) model. As reviewed by many authors (see, e.g., Yang et al. (2022, 2023b) and Li et al. (2023)), the standardized residual is an important indicator for an adequate model if the residual sequence $\{\hat{e}_{j,t}\}$ does not exhibit significant autocorrelation. Specifically, if the model is correctly specified, the residuals should have no significant serial correlation. The standardized residual plot can be used to examine the independent and identically distributed assumption of the residuals and to detect possible outliers for the dataset. To show the fitting details of the I-BRCMNBINAR(1) model, we draw the diagnostic checking plots in Fig. \ref{Pearson_residual} based on standardized Pearson residuals, which are defined by
\begin{align*}
e_{j,t}=\frac{X_{j,t}-E(X_{j,t}|X_{j,t-1})}{\sqrt{Var(X_{j,t}|X_{j,t-1})}},\ j=1,2.
\end{align*}
In practice, we calculate $\{\hat{e}_{j,t}\}$ by substituting the CML-estimator into the conditional expectation and conditional variance equations. From Fig. \ref{Pearson_residual} we can see that most lags of ACF and PACF values are within the blue dotted lines. The mean and variance for the two residual sequences are (0.00003, 0.05356) and (0.514,1.092). To verify the two residuals series are stationary white noise, we carry out the ADF test and Ljung-Box (LB) test for the fitted standard Pearson residuals. The $p$-values of the ADF test for the two residual series are both smaller than 0.01. We have $Q(12)=18.805(0.09)$ and $Q(12)=10.203(0.59)$ for the LB test of the standardized residuals of the two residual series respectively, where the number in the parentheses denotes $p$-value. The results show that $\{\hat{e}_{j,t}\}(j\in\{1,2\})$ are stationary white noise which ensure the I-BRCMNBINAR(1) model is correctly specified. We also consider introducing covariates in the innovation terms proposed in Yang et al. (2023b) to further describe the serial dependence of the data in a future research.
%However, this issue is beyond the scope of this paper at this stage. Therefore, we leave it as a topic of a future research.

\begin{figure}[!t]
\centering
\includegraphics[width=14cm,height=5cm]{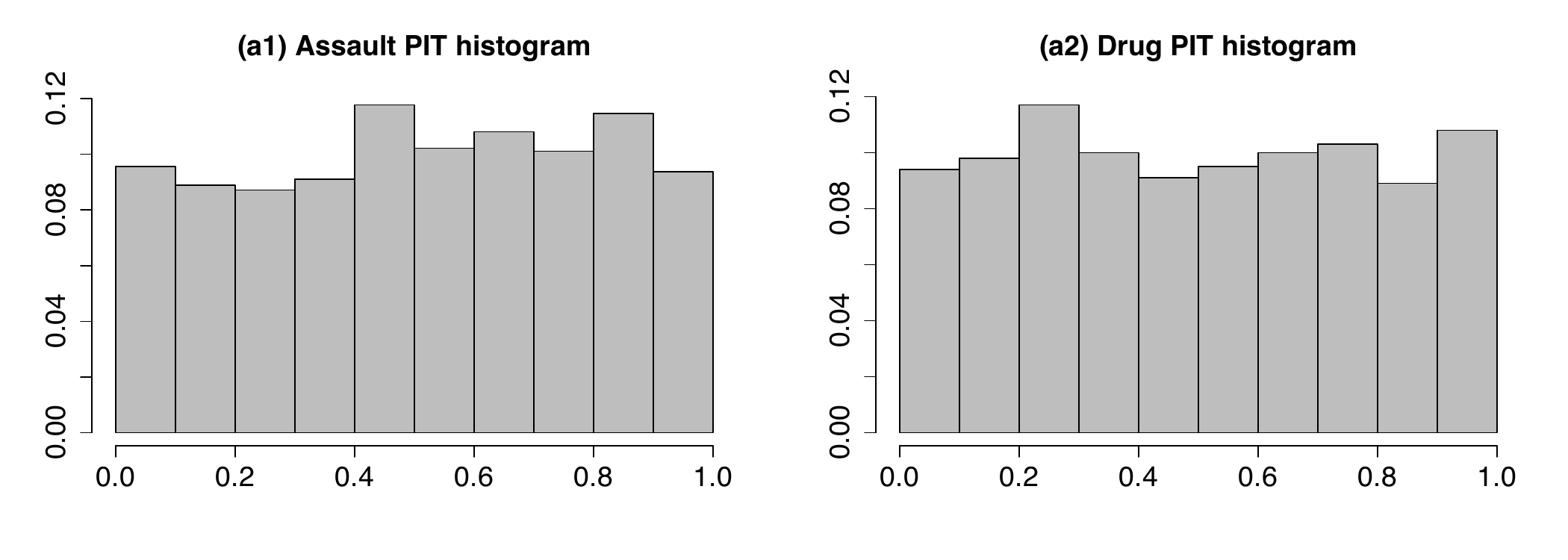}
\includegraphics[width=14cm,height=5cm]{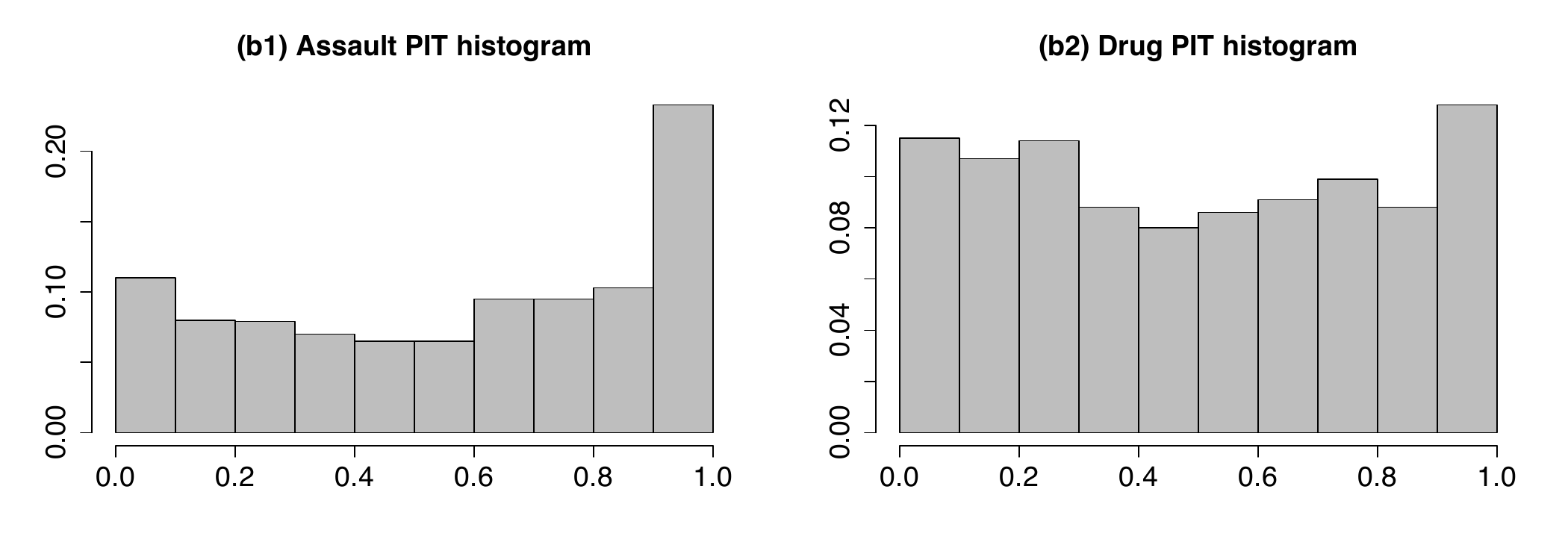}
\caption{PIT histograms of I-BRCMNBINAR(1) model with BVNB innovations (top figures) and BP innovations (bottom figures).}
\label{PIT}
\end{figure}

%p value???
To assess the adequacy of different distributional form for the innovations, we considered the probability integral transform (PIT) histogram proposed by Czado et al. (2009). It is widely used to compare the fitting effect of various models (see e.g. Jung et al. (2016), Monteiro et al. (2021) and Qian et al. (2020)). As argued by Wei\ss (2018), we expect the histogram of the series of PIT random variables to look like that of a uniform distribution if the fitted model is adequate. Fig. \ref{PIT} displays the nonrandomized PIT histograms. As we can see that the PIT histogram of the BVNB innovations is close to uniformity, which means it performs better than the BP innovations. In conclusion, based on the value of AIC, BIC, the analysis of the standard Pearson residuals and the nonrandomized PIT histograms, I-BRCMNBINAR(1) model shows the best performance. Therefore, these results suggest that it is more competitive to use the I-BRCMNBINAR(1) model to fit the dataset.
%competitive/superior/appropriate

%In PIT histograms approximate 95\% confidence intervals, obtained from a standard $\chi^2$ goodness-of-fit test (being the null hypothesis that the $J=10$ bins of the histogram are drawn from a uniform distribution) were incorporated as in Jung, McCabe, and Tremaayne (2015). 
%Under correct specification of the predictive distribution, the series of PIT random variables {ut} are i.i.d. standard uniform (0, 1). 
%As noted in [27], under correct specification of predictive distribution, the series of PIT random variables are i.i.d. standard uniform (0,1).
%Figure 5 also shows cumulative predictive distribution against the uniform distribution.
%The forecast RMSE and PIT histograms are shown to indicate that our proposed model is appropriate to model this data set and has a good performance in prediction and model diagnostics. 

\section{Possible extension}\label{SEC7}
Inspired by Monteiro et al. (2021) and Yang et al. (2023b), model (\ref{eq2}) can be extended to a new multivariate random coefficient threshold INAR process based on modified negative binomial operator with time-dependent innovation vectors, i.e.
\begin{align}\label{multi_threshold}
\textit{\textbf{X}}_{t}=(\textit{\textbf{A}}_t\diamond\textit{\textbf{X}}_{t-1})I_{1,t}(r)+(\textit{\textbf{B}}_t\diamond\textit{\textbf{X}}_{t-1})I_{2,t}(r)+\textit{\textbf{Z}}_{t},\; t\in\mathbb{N},
\end{align}
where $\textit{\textbf{X}}_{t}=(X_{1,t},\cdots,X_{K,t})^T\in \mathbb{N}_0^K$, $\textit{\textbf{A}}_t=diag\{\alpha_{1,t},\cdots,\alpha_{K,t}\}$ and $\textit{\textbf{B}}_t=diag\{\beta_{1,t},\cdots,\beta_{K,t}\}$. For $j\in\{1,2,\cdots,K\}$, the diagonal elements $\alpha_{j,t}$ and $\beta_{j,t}$ follow the Beta prime distribution with parameters $m_{1,j}$ and $m_{2,j}$ given that the known constant $l$. Let $R_t$ denote the threshold variable, $I_{1,t}(r)=I\{R_t\leqslant r\}$, $I_{2,t}(r)=I\{R_t> r\}$, $r$ is the unknown threshold parameter. The $\mathbb{N}_0^K$-valued innovation vector $\textit{\textbf{Z}}_{t}=(Z_{1,t},\cdots, Z_{K,t})^T$ follows one particular multivariate distribution. For fixed $t$ and $j$ ($j\in \{1,2,\cdots,K\}$), $\textit{\textbf{Z}}_{t}$ is assumed to be independent of $\textit{\textbf{A}}_s\diamond\textit{\textbf{X}}_{t-s}$ and $\textit{\textbf{B}}_s\diamond\textit{\textbf{X}}_{t-s}$ for $s<t$.

In fact, model (\ref{multi_threshold}) is a two-regime threshold multivariate autoregressive process. $\textit{\textbf{X}}_{t}$ takes the first regime if $R_t-r\leqslant 0$ and $\textit{\textbf{X}}_{t}$ takes the second regime if $R_t-r> 0$. For $i=1,2$, we omit $``(r)"$ in $``I_{i,t}(r)"$ to make the notations more clear. As suggested by Yang et al. (2023b), $R_t$ can be selected as a function of $\textit{\textbf{X}}_{t-1}$, such as $R_t=X_{1,t-1}$ or $R_t=\max\{X_{1,t-1},\cdots,X_{K,t-1}\}$. We assume that $\textit{\textbf{A}}_{t}\neq \textit{\textbf{B}}_{t}$ and $0<\underline{r}\leqslant r \leqslant\bar{r}$. Otherwise, the existence of the threshold variable is meaningless. The univariate process for each component of $\textit{\textbf{X}}_t$ takes the following form
\begin{align*}
X_{j,t}=\delta_{j,t}\diamond X_{j,t-1}+Z_{j,t},\ t\in\mathbb{N},
\end{align*}
where $\delta_{j,t}=\alpha_{j,t}I_{1,t}+\beta_{j,t}I_{2,t},\ j=1,2,\cdots,K$.

Now we need to consider the distribution of $\textit{\textbf{Z}}_{t}$. For example, Yang et al. (2023b) proposed a multivariate threshold integer-valued autoregressive process with explanatory variables, where the innovation vector follows a multivariate Poisson distribution. To capture time dependence, Chen et al. (2022) introduced an extended bivariate INAR(1) process where the mean of the innovation vector is linearly increased by the previous population size. As the mixed negative binomial models have thick tails, Tzougas and Cerchiara (2021) proposed to consider it to capture overdispersion and positive dependencies for multivariate count data, As an illustration, we consider a special case of this model, that is, 
\begin{align*}
Z_{j,t}\sim NB(\vartheta_{j},\frac{\vartheta_{j}}{\vartheta_{j}+\psi_{j,t}}),
\end{align*}
where $\psi_{j,t}=\exp(\textbf{\textit{X}}_{t-1}^T\pmb{\zeta}_j)$, $\pmb{\zeta}_j=(\zeta_j^{(1)},\cdots,\zeta_j^{(K)})^T$ contains the values of regression coefficients, $j=1,2,\cdots,K$. It is easy to verify that the process $\{\textit{\textbf{X}}_t\}_{t\in\mathbb{N}}$ generated by the MRCTMNBINAR(1) process defined by (\ref{multi_threshold}) is a Markov chain on state space $\mathbb{N}_0^K$ and the conditional expectation is given by:
\begin{align*}
&E(X_{j,t}|\textit{\textbf{X}}_{t-1})=\sum_{i=1}^2m_{i,j}(X_{j,t-1}+1)I_{i,t}/l+\psi_{j,t}.%\\
%&Var(X_{j,t}|\textit{\textbf{X}}_{t-1})= \\
%&Cov(X_{i,t},X_{j,t}|\textit{\textbf{X}}_{t-1})
\end{align*}
Furthermore, similar with Monteiro et al. (2021) and Chen et al. (2022), the transition probabilities of the MRCTMNBINAR(1) process are given by:
\begin{align*}
P(\textit{\textbf{X}}_{t}=\textit{\textbf{x}}_{t}|\textit{\textbf{X}}_{t-1}=\textit{\textbf{x}}_{t-1})
=\sum_{s_1=0}^{x_{1,t}}\cdots\sum_{s_K=0}^{x_{K,t}} 
\prod_{j=1}^{K}p_j(s_j)f_{\textit{\textbf{Z}}}(\textit{\textbf{x}}_{t}-\textit{\textbf{s}}),
\end{align*}
where $\textit{\textbf{s}}=(s_1,\cdots,s_K)^T$, $f_{\textit{\textbf{Z}}}(\textit{\textbf{x}}_{t}-\textit{\textbf{s}})$ denotes the probability mass function of $\textit{\textbf{Z}}_{t}$ and 
\begin{align*}
p_j(s_j)
=\sum_{i=1}^{2}\binom{x_{j,t-1}+s_j}{s_j}\frac{\Gamma(x_{j,t-1}+l+2)\Gamma(s_j+m_{i,j})\Gamma(l+m_{i,j}+1)I_{i,t}}{\Gamma(x_{j,t-1}+s_j+m_{i,j}+l+2)\Gamma(l+1)\Gamma(m_{i,j})}.
\end{align*}

The MRCTMNBINAR(1) model defined in (\ref{multi_threshold}) is a generalization of the BRCMNBINAR(1) model. It is able to capture the time-dependence trend by imposing the past information in the distribution of the innovation vector and introduce the cross-correlation between the multiple components into the innovation vector. The stationarity and ergodicity of process (\ref{multi_threshold}) need to be discussed, which is critical to derive the consistency and asymptotic normality of the CML estimator. However, these issues on the MRCTMNBINAR(1) model is somewhat beyond the scope of this paper at this stage and require further attention, we leave it as a future project.

\section{Conclusions}\label{SEC8}
%, especially for large sample sizes
This article proposes a new bivariate random coefficient INAR(1) process based on modified negative binomial operator with dependent innovations. The stationarity and ergodicity of the process are established. The YW-estimator, CLS-estimator and CML-estimator are derived and the related asymptotic properties are obtained. As an illustration, we conduct a simulation study to examine the effectiveness of these methods and the result shows that it is more reliable to use the CML method to estimate the parameters because of its great effect. The coherent forecasts for the BRCMNBINAR(1) model are also discussed. Monthly crime datasets in Bellingen are analyzed. The fitting result reveals that our proposed model can better describes the pattern of datasets, which illustrates the practicality of our model. Further generalization in the future research might concern the introduction of threshold variable and other distribution of the innovations in the multivariate time series.

\section*{Acknowledgement}
This work is supported by Social Science Planning Foundation of Liaoning Province (No. L22ZD065), National Natural Science Foundation of China (No. 12271231, 12001229, 11901053) and China Scholarship Council (Grant No. CSC202206170056).

\section*{Appendix}
{\bf Proof of Proposition \ref{proposition2.3}}.

We first introduce a bivariate random sequence $\{\textbf{\textit{X}}_{t}^{(n)}\}_{n\in\mathbb{N}_0}$ as follows:
\begin{align}\label{introduce_sequence}
\textbf{\textit{X}}_{t}^{(n)}
=\left\{
\begin{aligned}
&0,\qquad\qquad\qquad\qquad\qquad\quad n<0,\\
&\textbf{\textit{Z}}_{t},\qquad\qquad\qquad\qquad\qquad\; n=0,\\
&\textbf{\textit{A}}_{t}\diamond\textbf{\textit{X}}_{t-1}^{(n-1)}+\textbf{\textit{Z}}_{t},\qquad\qquad n>0.
\end{aligned}\right.
\end{align}

For $\forall n\in\mathbb{N}_0$, $\textbf{\textit{Z}}_{t}$ is independent of $\textbf{\textit{A}}_{t}\diamond\textbf{\textit{X}}_{t-1}^{(n-1)}$ and $\textbf{\textit{X}}_{s}^{(n)},\ s<t$. Let $L^2(\Omega,\mathcal{F},P)=\{\textbf{\textit{X}}|E(\textbf{\textit{XX}}^T)<\infty\}$ denote the Hilbert space, where the measure on $L^2(\Omega,\mathcal{F},P)$ is given by $d(\textbf{\textit{X}},\textbf{\textit{Y}})=E(\textbf{\textit{XY}}^T)$.\\

(A1) Existence.

\textbf{Step 1:} $\{\textbf{\textit{X}}_{t}^{(n)}\}_{n\in\mathbb{N}_0}\in L^2(\Omega,\mathcal{F},P)$.

We begin with the expectation of $\textbf{\textit{X}}_{t}^{(n)}$. Let $\pmb{\mit\mu}^{(n)}$ denote the expectation vector of $\textbf{\textit{X}}_{t}^{(n)}$, then we have
\begin{align*}
\pmb{\mit\mu}^{(n)}&=E(\textbf{\textit{X}}_{t}^{(n)})
=E(\textbf{\textit{A}}_{t}\diamond\textbf{\textit{X}}_{t-1}^{(n-1)}+\textbf{\textit{Z}}_{t})\\
&=\textbf{\textit{A}}E(\textbf{\textit{X}}_{t-1}^{(n-1)})+\textbf{\textit{A}}\textbf{\textit{e}}+\pmb{\mit\mu}_{Z}
=\sum_{i=1}^{n}\textbf{\textit{A}}^i\textbf{\textit{e}}+\sum_{i=0}^{n}\textbf{\textit{A}}^i\pmb{\mit\mu}_{Z}, 
\end{align*}
where $\pmb{\mit\mu}_{Z}=(\mu_1,\mu_2)^T,\ \textbf{\textit{e}}=(1,1)^T$. Since $\lim_{n\to\infty}\sum_{i=0}^{n}\textbf{\textit{A}}^i=(\textbf{\textit{I}}-\textbf{\textit{A}})^{-1}$ and $\det|\textbf{\textit{I}}-\textbf{\textit{A}}|=(1-\frac{m_1}{l})(1-\frac{m_2}{l})>0$, then we have 
\begin{align*}
\pmb{\mit\mu}^{(n)}=(\textbf{\textit{I}}-\textbf{\textit{A}})^{-1}(\textbf{\textit{A}}\textbf{\textit{e}}+\pmb{\mit\mu}_{Z})<\infty,\ as\ n\rightarrow\infty.
\end{align*}

%要证二阶矩<\infty
For the second moments of $\textbf{\textit{X}}_{t}^{(n)}$, from Proposition \ref{proposition2.2}, we have 
\begin{align}%注释原因：排版美观
E[\textbf{\textit{X}}_{t}^{(n)}(\textbf{\textit{X}}_{t}^{(n)})^T]%\nonumber\\
%&=E[(\textbf{\textit{A}}_t\diamond\textbf{\textit{X}}_{t-1}^{(n-1)}+\textbf{\textit{Z}}_t)(\textbf{\textit{A}}_t\diamond\textbf{\textit{X}}_{t-1}^{(n-1)}+\textbf{\textit{Z}}_t)']\nonumber\\
%&=E[(\textbf{\textit{A}}_t\diamond\textbf{\textit{X}}_{t-1}^{(n-1)})(\textbf{\textit{A}}_t\diamond\textbf{\textit{X}}_{t-1}^{(n-1)})']+E[(\textbf{\textit{A}}_t\diamond\textbf{\textit{X}}_{t-1}^{(n-1)})\textbf{\textit{Z}}_t']+E[\textbf{\textit{Z}}_t(\textbf{\textit{A}}_t\diamond\textbf{\textit{X}}_{t-1}^{(n-1)})']+E(\textbf{\textit{Z}}_t\textbf{\textit{Z}}_t')\nonumber\\
%&=\textbf{\textit{A}}E[(\textbf{\textit{X}}_{t-1}^{(n-1)}+\textbf{\textit{e}})(\textbf{\textit{X}}_{t-1}^{(n-1)}+\textbf{\textit{e}})']\textbf{\textit{A}}+\textbf{\textit{C}}+\textbf{\textit{A}}(E(\textbf{\textit{X}}_{t-1}^{(n-1)})+\textbf{\textit{e}})\pmb{\mu}_{Z}'\nonumber\\
%&\qquad+\pmb{\mu}_{Z}(\textbf{\textit{e}}'+E(\textbf{\textit{X}}_{t-1}^{(n-1)})')\textbf{\textit{A}}+E(\textbf{\textit{Z}}_t\textbf{\textit{Z}}_t')\nonumber\\
&=\textbf{\textit{A}}E[\textbf{\textit{X}}_{t-1}^{(n-1)}(\textbf{\textit{X}}_{t-1}^{(n-1)})^T]\textbf{\textit{A}}
+\textbf{\textit{A}}E(\textbf{\textit{X}}_{t-1}^{(n-1)})\textbf{\textit{e}}^T\textbf{\textit{A}}+\textbf{\textit{A}}\textbf{\textit{e}}E(\textbf{\textit{X}}_{t-1}^{(n-1)})^T\textbf{\textit{A}}\nonumber\\
&\qquad+\textbf{\textit{A}}\textbf{\textit{e}}\textbf{\textit{e}}^T\textbf{\textit{A}}+\textbf{\textit{C}}+\textbf{\textit{A}}E(\textbf{\textit{X}}_{t-1}^{(n-1)})\pmb{\mu}_{Z}^T+\textbf{\textit{A}}\textbf{\textit{e}}\pmb{\mu}_{Z}^T\nonumber\\
&\qquad+\pmb{\mu}_{Z}\textbf{\textit{e}}^T\textbf{\textit{A}}+\pmb{\mu}_{Z}E(\textbf{\textit{X}}_{t-1}^{(n-1)})^T\textbf{\textit{A}}+E(\textbf{\textit{Z}}_t\textbf{\textit{Z}}_t^T).
\label{eq5}
\end{align}
Notice that (\ref{eq5}) is independent of $t$, iterating it $n$ times, then we have 
\begin{align}
E[\textbf{\textit{X}}_{t}^{(n)}(\textbf{\textit{X}}_{t}^{(n)})^T]=\textbf{\textit{A}}^{n}E(\textbf{\textit{Z}}_t\textbf{\textit{Z}}_t^T)\textbf{\textit{A}}^n+\sum_{i=0}^{n-1}\textbf{\textit{A}}^{i}\textbf{\textit{M}}\textbf{\textit{A}}^{i},
\end{align}
where 
\begin{align*}
\textbf{\textit{M}}&=\textbf{\textit{A}}(\textbf{\textit{I}}-\textbf{\textit{A}})^{-1}(\textbf{\textit{Ae}}+\pmb{\mit\mu}_{Z})(\textbf{\textit{e}}^T\textbf{\textit{A}}+\pmb{\mu}_{Z}^T)+(\textbf{\textit{Ae}}+\pmb{\mu}_{Z})(\textbf{\textit{e}}^T\textbf{\textit{A}}+\pmb{\mit\mu}_{Z}^T)[(\textbf{\textit{I}}-\textbf{\textit{A}})^{-1}]^T\textbf{\textit{A}}\nonumber\\
&+\textbf{\textit{Aee}}^T\textbf{\textit{A}}+\textbf{\textit{C}}
+\textbf{\textit{Ae}}\pmb{\mu}_{Z}^T+\pmb{\mu}_{Z}\textbf{\textit{e}}^T\textbf{\textit{A}}+E(\textbf{\textit{Z}}_t\textbf{\textit{Z}}_t^T).
\end{align*}
As $\pmb{\mit\mu}^{(n)}=E(\textbf{\textit{X}}_{t}^{(n)})=(\textbf{\textit{I}}-\textbf{\textit{A}})^{-1}(\textbf{\textit{A}}\textbf{\textit{e}}+\pmb{\mit\mu}_{Z})<\infty$, we can conclude that $E[\textbf{\textit{X}}_{t}^{(n)}(\textbf{\textit{X}}_{t}^{(n)})^T]<\infty$ which implies $\{\textbf{\textit{X}}_{t}^{(n)}\}_{n\in\mathbb{N}_0}\in L^2(\Omega,\mathcal{F},P)$.\\

\textbf{Step 2:} $\{\textbf{\textit{X}}_{t}^{(n)}\}_{n\in\mathbb{N}_0}$ is a Cauchy sequence.

%要证非减
For $n=1$, we have 
\begin{align*}
\textbf{\textit{X}}_{t}^{(1)}
=\textbf{\textit{A}}_t\diamond\textbf{\textit{X}}_{t-1}^{(0)}+\textbf{\textit{Z}}_{t}
=\textbf{\textit{A}}_t\diamond\textbf{\textit{Z}}_{t}+\textbf{\textit{Z}}_{t}\geqslant\textbf{\textit{Z}}_{t}=\textbf{\textit{X}}_{t}^{(0)}.
\end{align*}
For $\forall t\in\mathbb{N},k=1,2,3,\cdots,n$, suppose that $\textbf{\textit{X}}_t^{(k)}\geqslant\textbf{\textit{X}}_t^{(k-1)}$. The $i$th element of $\textbf{\textit{X}}_{t}^{(n)}-\textbf{\textit{X}}_{t}^{(n-1)}$ for $n\geqslant 1$ is given by
\begin{align*}
\big(\textbf{\textit{X}}_{t}^{(n)}-\textbf{\textit{X}}_{t}^{(n-1)}\big)_i=\alpha_{i,t}\diamond X_{i,t-1}^{(n-1)}-\alpha_{i,t}\diamond X_{i,t-1}^{(n-2)}=\sum_{j=1}^{X_{i,t-1}^{(n-1)}-X_{i,t-1}^{(n-2)}}G_{ij}^{(t)}\geqslant 0,\; i=1,2.
\end{align*}
Therefore, $\textbf{\textit{X}}_{t}^{(n)}-\textbf{\textit{X}}_{t}^{(n-1)}\geqslant 0$ holds for $\forall t\in\mathbb{N}$, which implies $\{\textbf{\textit{X}}_t^{(n)}\}_{n\in\mathbb{N}_0}$ is a non-decreasing sequence.

Let $\textbf{\textit{U}}(n,t,k)=\textbf{\textit{X}}_t^{(n)}-\textbf{\textit{X}}_t^{(n-k)},\ k=1,2,\cdots,n$. Then we have
\begin{align*}
\textbf{\textit{U}}(n,t,k)
&=\left(U_1(n,t,k), U_2(n,t,k)\right)^T\\
&=\left(\sum_{j=1}^{X_{1,t-1}^{(n-1)}-X_{1,t-1}^{(n-k-1)}}G_{1j}^{(t)},\sum_{j=1}^{X_{2,t-1}^{(n-1)}-X_{2,t-1}^{(n-k-1)}}G_{2j}^{(t)}\right)^T\\
&=\left(\sum_{j=1}^{U_1(n-1,t-1,k)}G_{1j}^{(t)},\sum_{j=1}^{U_2(n-1,t-1,k)}G_{2j}^{(t)}\right)^T.
\end{align*}

%很重要: EU(n,t,k)=0!!!
For the expectation of $\textbf{\textit{U}}(n,t,k)$, we have 
\begin{align*}
E[\textbf{\textit{U}}(n,t,k)]=\textbf{\textit{A}}E[\textbf{\textit{U}}(n-1,t-1,k)]=\textbf{\textit{A}}^n\pmb{\mit\mu}_{Z}\rightarrow\textbf{\textit{0}},\; n\rightarrow\infty.
\end{align*}
For $i=1,2$, we have
\begin{align*} 
E[U_i(n,t,k)^2]
&=E\left[\left(\sum_{j=1}^{U_i(n-1,t-1,k)}G_{ij}^{(t)}\right)^2\right]\\  
&=(\sigma^2_{\alpha_i}+\alpha_i^2)E[U_i(n-1,t-1,k)^2]+(\sigma^2_{\alpha_i}+\alpha_i^2+\alpha_i)E[U_i(n-1,t-1,k)]\\
&=(\sigma^2_{\alpha_i}+\alpha_i^2+\alpha_i)\big((\sigma^2_{\alpha_i}+\alpha_i^2)E[U_i(n-2,t-2,k)]+E[U_i(n-1,t-1,k)]\big)\\
&\quad+(\sigma^2_{\alpha_i}+\alpha_i^2)^2E[U_i(n-2,t-2,k)^2]\\
&\cdots\\
&=(\sigma^2_{\alpha_i}+\alpha_i^2)^nE[U_i(0,t-n,k)^2]+(\sigma^2_{\alpha_i}+\alpha_i^2+\alpha_i)\frac{1-(\sigma^2_{\alpha_i}+\alpha_i^2)^n}{1-\sigma^2_{\alpha_i}-\alpha_i^2}\alpha_i^n\mu_{i}. 
\end{align*} 
Thus it is easy to see that $E[U_i(n,t,k)^2]\rightarrow 0$, as $n\rightarrow\infty$. Along the same recursive line, we have
\begin{align*} 
E[U_1(n,t,k)U_2(n,t,k)]=\alpha_1^n\alpha_2^n(\phi+\mu_{1}\mu_{2})\rightarrow 0,\; as\;n\rightarrow\infty.
\end{align*}
Therefore, we can conclude that $E[\textbf{\textit{U}}(n,t,k)\textbf{\textit{U}}(n,t,k)^T]\rightarrow 0,\; as\; n\rightarrow\infty$, which means $\{\textbf{\textit{X}}_{t}^{(n)}\}_{n\in\mathbb{N}_0}$ is a Cauchy sequence.\\

(A2) Uniqueness.

For $i=1,2$, suppose that there is another process $\{X_{i,t}^*\}_{t\in\mathbb{N}}$ such that $X_{i,t}^{(n)}\overset{L^2}{\rightarrow}X_{i,t}^*$. By the H\"{o}lder inequality, we have
\begin{align*}
E|X_{i,t}-X_{i,t}^*|\leqslant\big[E|X_{i,t}-X_{i,t}^{(n)}|^2\big]^{1/2}\big[E|X_{i,t}^*-X_{i,t}^{(n)}|^2\big]^{1/2}\rightarrow 0,\;as\; n\rightarrow\infty.
\end{align*}
Thus we have $E|X_{i,t}-X_{i,t}^*|=0$, which implies $X_{i,t}=X_{i,t}^*$ almost surely.\\

(A3) Strict stationarity.

Notice that $\{\textit{\textbf{A}}_t\}$ and $\{\textit{\textbf{Z}}_t\}$ are independent and identically distributed random sequences, and repeat application of equation (\ref{introduce_sequence}) with $n$ times gives
\begin{align*} 
\textit{\textbf{X}}_t^{(n)}
&=\textit{\textbf{A}}_t\diamond\textit{\textbf{A}}_{t-1}\cdots\textit{\textbf{A}}_{1}\diamond\textbf{\textit{Z}}_{t-n}+\sum_{i=1}^{n-1}\textit{\textbf{A}}_{t}\diamond\cdots\diamond\textit{\textbf{A}}_{t-i+1}\diamond\textbf{\textit{Z}}_{t-i}+\textbf{\textit{Z}}_{t}\\
%用了同分布, 把t都换成n, see zhang qingchun2020
&\overset{d}{=}\textit{\textbf{A}}_n\diamond\textit{\textbf{A}}_{n-1}\cdots\textit{\textbf{A}}_{1}\diamond\textbf{\textit{Z}}_{0}
+\sum_{i=1}^{n-1}\textit{\textbf{A}}_{n}\diamond\cdots\diamond\textit{\textbf{A}}_{n-i+1}\diamond\textbf{\textit{Z}}_{n-i}+\textbf{\textit{Z}}_{n},
\end{align*}    
where $\overset{d}{=}$ denotes equality in distribution. Thus, the distribution of $\textbf{\textit{X}}_t^{(n)}$ depends only on $n$ rather than $t$. For $\forall k,T\in\mathbb{N}$, we have 
\begin{align*}
(\textbf{\textit{X}}_1^{(n)},\cdots,\textbf{\textit{X}}_T^{(n)})\overset{d}{=}(\textbf{\textit{X}}_k^{(n)},\cdots,\textbf{\textit{X}}_{k+T-1}^{(n)}). 
\end{align*} 

Apply the Cram\'{e}r-Wold Theorem, we have
\begin{align*}
(\textbf{\textit{X}}_1^{(n)},\cdots,\textbf{\textit{X}}_T^{(n)})\xrightarrow{L}(\textbf{\textit{X}}_1,\cdots,\textbf{\textit{X}}_T),\ (\textbf{\textit{X}}_k^{(n)},\cdots,\textbf{\textit{X}}_{k+T-1}^{(n)})\xrightarrow{L}(\textbf{\textit{X}}_k,\cdots,\textbf{\textit{X}}_{k+T-1}).
\end{align*} 
Therefore, we can conclude that $(\textbf{\textit{X}}_1,\cdots,\textbf{\textit{X}}_T)\overset{d}{=}(\textbf{\textit{X}}_k,\cdots,\textbf{\textit{X}}_{k+T-1})$, which means $\{\textbf{\textit{X}}_t^{(n)}\}_{n\in\mathbb{N}_0}$ is a strictly stationary process. As the sequence $\{\textbf{\textit{X}}_t^{(n)}\}_{n\in\mathbb{N}_0}$ converges to $\{\textbf{\textit{X}}_t\}_{t\in\mathbb{N}}$ in $L^2(\Omega,\mathcal{F},P)$, the process $\{\textbf{\textit{X}}_t\}_{t\in\mathbb{N}}$ is a strictly stationary process.\\

(A4) Ergodicity.

Let $\sigma(\textbf{\textit{X}})$ be the $\sigma$-field generated by the random vector $\textbf{\textit{X}}$ and $\textbf{\textit{W}}_t$ be all counting series involved in random matricial operation $\textbf{\textit{A}}_t\diamond\textbf{\textit{X}}_{t-1}$. It is obviously that $\textbf{\textit{W}}_t$ is a time series. From (\ref{eq2}), we have 
\begin{align*}
\sigma(\textbf{\textit{X}}_t,\textbf{\textit{X}}_{t+1},\cdots)\subset\sigma(\textbf{\textit{Z}}_t,\textbf{\textit{A}}_t,\textbf{\textit{W}}_t,\textbf{\textit{Z}}_{t+1},\textbf{\textit{A}}_{t+1},\textbf{\textit{W}}_{t+1},\cdots).
\end{align*} 

Thus, we can obtain that
\begin{align*}
\bigcap_{t=1}^{\infty}\sigma(\textbf{\textit{X}}_t,\textbf{\textit{X}}_{t+1},\cdots)\subset\bigcap_{t=1}^{\infty}\sigma(\textbf{\textit{Z}}_t,\textbf{\textit{A}}_t,\textbf{\textit{W}}_t,\textbf{\textit{Z}}_{t+1},\textbf{\textit{A}}_{t+1},\textbf{\textit{W}}_{t+1},\cdots).
\end{align*}  
Since $\{(\textbf{\textit{Z}}_t,\textbf{\textit{A}}_t,\textbf{\textit{W}}_t)\}$ is an independent random vectors sequence, Kolmogorov's zero-one law implies that for any event $A\in\bigcap_{t=1}^{\infty}\sigma(\textbf{\textit{Z}}_t,\textbf{\textit{A}}_t,\textbf{\textit{W}}_t,\textbf{\textit{Z}}_{t+1},\textbf{\textit{A}}_{t+1},\textbf{\textit{W}}_{t+1},\cdots)$, we have $P(A)=0$ or $P(A)=1$. Thus $\{\textbf{\textit{X}}_t\}_{t\in\mathbb{N}}$ is an ergodic process as the tail of the $\sigma$-field of $\{\textbf{\textit{X}}_t\}_{t\in\mathbb{N}}$ contains only the measure sets with probability 0 or 1 known from Wang (1982).\\

{\bf Proof of Proposition \ref{proposition2.4}}.

The expectations are easy to verfied. We begin with the conditional variance of $X_{i,t+1}$ on $X_{1,t}$ and $X_{2,t}$,  
\begin{align*}
V&ar(X_{i,t+1}|X_{1,t},X_{2,t})\\
&=E_{\alpha_{i,t+1}}\big(Var(\alpha_{i,t+1}\diamond X_{i,t}+Z_{i,t+1}|X_{i,t},\alpha_{i,t+1})\big)\\
&\quad+Var_{\alpha_{i,t+1}}\big(E(\alpha_{i,t+1}\diamond X_{i,t}+Z_{i,t+1}|X_{i,t},\alpha_{i,t+1})\big)\\
&=E_{\alpha_{i,t+1}}\big((X_{i,t}+1)\alpha_{i,t+1}(1+\alpha_{i,t+1})\big)+\sigma^2_{Z_i}+Var_{\alpha_{i,t+1}}\big((X_{i,t}+1)\alpha_{i,t+1}\big)\\
&=(\sigma^2_{\alpha_i}+\alpha_i^2+\alpha_i)(X_{i,t}+1)+\sigma^2_{Z_i}+(X_{i,t}+1)^2\sigma^2_{\alpha_i}.
\end{align*}
From law of total variance, we have
\begin{align*}
Var(X_{i,t})&=Var\big(E(X_{i,t+1}|X_{i,t})\big)+E\big(Var(X_{i,t+1}|X_{i,t})\big)\\
&=\alpha_i^2Var(X_{i,t}+1)+\sigma^2_{Z_i}+(\sigma^2_{\alpha_i}+\alpha_i^2+\alpha_i)E(X_{i,t}+1)+\sigma^2_{\alpha_i}E(X_{i,t+1}+1)^2.
\end{align*} 
Then we can derive the unconditional variance of $X_{i,t}$.\\

For one hand, the covariance of $X_{1,t}$ and $X_{2,t}$ is given by
\begin{align*}
Cov(X_{1,t},X_{2,t})&=Cov(\alpha_{1,t}\diamond X_{1,t-1},\alpha_{2,t}\diamond X_{2,t-1})+Cov(Z_{1,t},Z_{2,t})\\
&=\frac{m_1m_2}{l^2}Cov(X_{1,t-1},X_{2,t-1})+Cov(Z_{1,t},Z_{2,t})\\
&=\frac{m_1^2m_2^2}{l^4}Cov(X_{1,t-2},X_{2,t-2})+\frac{m_1m_2}{l^2}Cov(Z_{1,t-1},Z_{2,t-1})+Cov(Z_{1,t},Z_{2,t})\\
&\cdots\\
&=\frac{l^2\phi}{l^2-m_1m_2}.
\end{align*}
For another hand, the covariance of $X_{i,t+k}$ and $X_{j,t}$ is given by
\begin{align*}
Cov(&X_{i,t+k},X_{j,t})\\
&=Cov(\alpha_{i,t+k}\diamond\cdots\alpha_{i,t+1}\diamond X_{i,t}+\sum_{n=1}^{k-1}\alpha_{i,t+k}\diamond\cdots\diamond\alpha_{i,t+k-n+1}\diamond Z_{i,t+k-n}+Z_{i,t+k},X_{j,t})\\
&=E((\alpha_{i,t+k}\diamond\cdots\alpha_{i,t+1}\diamond X_{i,t})X_{j,t})-E(\alpha_{i,t+k}\diamond\cdots\alpha_{i,t+1}\diamond X_{i,t})E(X_{j,t})\\
&=\alpha_i^kE(X_{i,t}X_{j,t})-\alpha_i^kE(X_{i,t})E(X_{j,t})\\
&=m_i^k\phi/(l^{k}-l^{k-2}m_1m_2).
\end{align*} 
This completes all the proof.\\

{\bf Proof of Theorem \ref{th3.1}}. 

The strong consistency of the CLS-estimator $\hat{\pmb{\gamma}}^{CLS}$ can be proved by checking the regularity conditions in Klimko and Nelson (1978) are satisfied, since the process $\{\textit{\textbf{X}}_t\}_{t\in\mathbb{N}}$ is stationary and ergodic. Now we consider the asymptotic normality. Let $\mathcal{F}_{n-1}=\sigma(\textit{\textbf{X}}_1,\cdots,\textit{\textbf{X}}_{n-1})$ denote the information filtration until time $n-1$. For $\textbf{\textit{D}}_n=\sum_{t=2}^{n}\frac{\partial g(\pmb{\gamma},{\textbf{\textit{X}}}_{t-1})^T}{\partial\pmb{\gamma}}\textbf{\textit{u}}_t$. Let $D_{11}=D_{12}=D_{13}=D_{14}=0$. For $n\geqslant 2$, $i=1,2$, we have 
\begin{align*}
E(D_{ni}|\mathcal{F}_{n-1})&=E\big((X_{i,n}-m_i(X_{i,n-1}+1)/l-\mu_i)(X_{i,n-1}+1)/l+D_{(n-1)i}|\mathcal{F}_{n-1}\big)\\
&=D_{(n-1)i}.
\end{align*}
Similarly, for $i=3,4$, we have 
\begin{align*}
E(D_{ni}|\mathcal{F}_{n-1})&=E\big((X_{i-2,n}-m_{i-2}(X_{i-2,n-1}+1)/l-\mu_{i-2})+D_{(n-1)i}|\mathcal{F}_{n-1}\big)=D_{(n-1)i}.
\end{align*}
Thus, $\{D_{ni},\mathcal{F}_{n},n\geqslant 1\}$ is a martingale. Furthermore, for any nonzero vector $\textbf{\textit{K}}=(K_1,K_2,K_3,K_4)^T\in\mathbb{R}^{4}\backslash(0,0,0,0)^T$, $\{\textbf{\textit{K}}^T\textbf{\textit{D}}_n,\mathcal{F}_{n},n\geqslant 1\}$ is a martingale. Since $E|X_{i,t}|^4<\infty,\ i=1,2$, $\{\textit{\textbf{X}}_t\}_{t\in\mathbb{N}}$ is a strictly stationary and ergodic process, we have 
\begin{align*}
\frac{1}{n-1}\sum_{t=2}^{n}\textbf{\textit{D}}_t\textbf{\textit{D}}_t^T\xrightarrow{a.s.}\textit{\textbf{W}}_1,\qquad\frac{1}{n-1}\sum_{t=2}^{n}(\textbf{\textit{K}}^T\textbf{\textit{D}}_t)^2\xrightarrow{a.s.}\sigma_1^2,
\end{align*}
where $\sigma_1^2=E(\textbf{\textit{K}}^T\textbf{\textit{D}}_n)^2$ and 
\begin{align} 
\textit{\textbf{W}}_1
=E\left(\frac{\partial g(\pmb{\gamma},\textbf{\textit{X}}_{t-1})^T\textbf{\textit{u}}_t\textbf{\textit{u}}_t^T\partial g(\pmb{\gamma},\textbf{\textit{X}}_{t-1})}{\partial{\pmb{\gamma}}\partial{\pmb{\gamma}}^T}\right)<\infty.
\end{align} 
Applying the martingale central limit theorem from Corollary 3.2 in Hall and Heyde (1980), as $n\rightarrow\infty$, we have 
\begin{align*}
\frac{1}{\sqrt{n-1}}\textbf{\textit{K}}^T\textbf{\textit{D}}_n
&=\frac{1}{\sqrt{n-1}}\sum_{t=2}^{n}(X_{1,t}-\alpha_1(X_{1,t-1}+1)-\mu_1)(K_1(X_{1,t-1}+1)/l+K_2)\\
&\qquad + (X_{2,t}-\alpha_2(X_{2,t-1}+1)-\mu_2)(K_3(X_{2,t-1}+1)/l+K_4)\xrightarrow{L}N(0,\sigma_1^2).
\end{align*}
According to the Cram\'{e}r-Wold device, we obtain 
\begin{align*}
\frac{1}{\sqrt{n-1}}\textbf{\textit{D}}_n=\frac{1}{\sqrt{n-1}}\sum_{t=2}^{n}\frac{\partial g(\pmb{\gamma},\textbf{\textit{X}}_{t-1})^T}{\partial\pmb{\gamma}}\textbf{\textit{u}}_t
\xrightarrow{L} 
N(\textit{\textbf{0}},\textit{\textbf{W}}_1).
\end{align*}
Since $\hat{\pmb{\gamma}}^{CLS}=\textit{\textbf{G}}_{n-1}^{-1}\textit{\textbf{b}}_{n-1},\ \textbf{\textit{V}}_1=\lim_{n\to\infty}\frac{1}{n-1}\textit{\textbf{G}}_{n-1}$, we have
\begin{align*}
\sqrt{n-1}(\hat{\pmb{\gamma}}^{CLS}-\pmb{\gamma})
=
\left(\frac{1}{n-1}\textit{\textbf{G}}_{n-1}\right)^{-1}
\frac{1}{\sqrt{n-1}}\textbf{\textit{D}}_n
\xrightarrow{L}
N(\textit{\textbf{0}},{\textit{\textbf{V}}_1}^{-1}{\textit{\textbf{W}}_1}{\textit{\textbf{V}}_1}^{-1}).
\end{align*}

{\bf Proof of Theorem \ref{th3.3}}.\\

Let $\textbf{\textit{B}}_n=(B_{n1},B_{n2},B_{n3},B_{n4})^T$ and $B_{11}=B_{12}=B_{13}=B_{14}=0$, then we have
\begin{align*}
B_{ni}&=-\frac{1}{2}\frac{\partial S_3(\pmb\eta)}{\partial\pmb{\eta}_i}=\sum_{t=2}^{n}U_{i,t}a_{i,t},\ n\geqslant 2, 
\end{align*}
with $U_{i,t}=(X_{i,t}-\hat{\alpha}_i(X_{i,t-1}+1)-\hat{\mu}_i)^2-\sigma^2_{\alpha_i}((X_{i,t-1}+1)^2+(X_{i,t-1}+1))-(\hat{\alpha}_i^2+\hat{\alpha}_i)(X_{i,t-1}+1)-\sigma^2_{Z_i}$ and $a_{i,t}=(X_{i,t}+1)^2+X_{i,t}+1,\ i=1,2$, and $a_{i,t}=1,\ i=3,4$. For $i=1,2,3,4$, we have 
\begin{align*}
E(B_{ni}|\mathcal{F}_{n-1})&=B_{(n-1)i}+E(U_{i,n}a_{i,n}|\mathcal{F}_{n-1})\\
&=B_{(n-1)i}+a_{i,n}E\big((X_{i,n}-E(X_{i,n}|X_{i,n-1}))^2-Var(X_{i,n}|X_{i,n-1})|\mathcal{F}_{n-1}\big)\\
&=B_{(n-1)i}.
\end{align*}
Therefore, $\{B_{ni},\mathcal{F}_{n},n\geqslant 1\}$ is a martingale. Using the ergodic theorem, we have
\begin{align*}
\frac{1}{n-1}\sum_{t=2}^{n}(U_{i,t}a_{i,t})^2\xrightarrow{a.s.} b_{ii},\\
\frac{1}{\sqrt{n-1}}B_{ni}\xrightarrow{L} N(0,b_{ii}),
\end{align*}
with $b_{ii}=E(U_{i,t}a_{i,t})^2<\infty$. Therefore, using the same arguments in the proof of Theorem \ref{th3.1} with $\textit{\textbf{V}}_2=E(\textbf{\textit{a}}_t\textbf{\textit{a}}_t^T)$ and $\textit{\textbf{W}}_2=E(\textbf{\textit{U}}_t^T\textbf{\textit{U}}_t\textbf{\textit{a}}_t\textbf{\textit{a}}_t^T)$ yields our desired result.


\begin{thebibliography}{99}
\bibitem[Aleksi\'{c} and Risti\'{c}(2021)]{Aleksić2021}
Aleksi\'{c}, M. S., and M. M. Risti\'{c}. 2021. A geometric minification integer-valued autoregressive model. {\textit{Applied Mathematical Modelling}} 90:265-280. doi:10.1016/j.apm.2020.08.047.
\bibitem[Al-Osh and Alzaid(1987)]{Al-Osh1987}%check
Al-Osh, M. A., and A. A. Alzaid. 1987. First-order integer-valued autoregressive (INAR(1)) process. {\textit{Journal of Time Series Analysis}} 8(3):261-275. doi:10.1111/j.1467-9892.1987.tb00438.x.


\bibitem[Billingsley(1961)]{Billingsley1961}%4
Billingsley, P. 1961. \textit{Statistical Inference for Markov Processes}. Chicago: The University of Chicago Press.


\bibitem[Chen et al(2022)]{Chen2022}
Chen, H., F. Zhu, and X. Liu. 2022. A new bivariate INAR(1) model with time-dependent innovation vectors. {\textit{Stats}} 5(3):819-840. doi:10.3390/stats5030048.
\bibitem[Czado et al(2009)]{Czado2009}
Czado, C., T. Gneiting, and L. Held. 2009. Predictive model assessment for count data. \textit{Biometrics} 65(4):1254-1261. doi:10.1111/j.1541-0420.2009.01191.x.



\bibitem[Ding and Wang(2016)]{Ding2016}%check
Ding, X., and D. Wang. 2016. Empirical likelihood inference for INAR(1) model with explanatory variables. {\textit{Journal of The Korean Statistical Society}} 45(4):623-632. doi:10.1016/j.jkss.2016.05.004.


\bibitem[Freeland and McCabe(2004)]{Freeland2004}%7
Freeland, R. K., and B. P. M. McCabe. 2004. Forecasting discrete-valued low count time series. {\textit{International Journal of Forecasting}} 20(4):427-434. doi:10.1016/S0169-2070(03)00014-1. 
 


\bibitem{Gorgi2020}
Gorgi, P. 2020. Beta-negative binomial auto-regressions for modelling integer-valued time series with extreme observations. {\textit{Journal of the Royal Statistical Society Series B: Statistical Methodology}} 82(5):1325-1347. doi:10.1111/rssb.12394.


\bibitem[Hall and Heyde(1980)]{Hall1980}%8
Hall, P., and C. C. Heyde. 1980. {\textit{Martingale Limit Theory and its Application}}. New York: Academic Press.
\bibitem[Hwang and Basawa(1998)]{Hwang1998}%9
Hwang, S. Y., and I. V. Basawa. 1998. Parameter estimation for generalized random coefficient autoregressive processes. {\textit{Journal of Statistical Planning and Inference}} 68(2):323-337. doi:10.1016/S0378-3758(97)00147-X.



\bibitem[Jung(2016)]{Jung2016}
Jung, R., B. McCabe, and A. Tremaayne. 2016. Model validation and diagnostics. In \textit{Handbook of
discrete-valued time series}, eds. R. A. Davis, H. Holan, R. Lund, and N. Ravishanker, 189-218. Boca Raton, FL: Chapman \& Hall.




\bibitem[Karlis and Pedeli(2013)]{Karlis2013}
Karlis, D., and X. Pedeli. 2013. Flexible bivariate INAR(1) processes using copulas. {\textit{Communications in Statistics-Theory and Methods}} 42(4):723-740. doi:10.1080/03610926.2012.754466.
\bibitem[Karlsen and Tj\o{}stheim(1988)]{Karlsen1988}
Karlsen, H., and D. Tj\o{}stheim. 1988. Consistent estimates for the NEAR(2) and NLAR(2) time series models. {\textit{Journal of the Royal Statistical Society Series B: Statistical Methodology}} 50(2):313-320. doi:10.1111/j.2517-6161.1988.tb01730.x.
\bibitem[Klenke(2013)]{Klenke2013}
Klenke, A. 2013. {\textit{Probability Theory: A Comprehensive Course}}. Springer Science \& Business Media.
\bibitem[Klimko and Nelson(1978)]{Klimko1978}
Klimko, L. A., and P. I. Nelson. 1978. On conditional least squares estimation for stochastic processes. {\textit{The Annals of Probability}} 6:629-642. doi:10.1214/aos/1176344207.



\bibitem{Li et al.(2023)}%check
Li, H., Z. Liu, K. Yang, X. Dong, and W. Wang. 2023. A $p$th-order random coefficients mixed binomial autoregressive process with explanatory variables. {\textit{Computational Statistics}} 1-24. doi:10.1007/s00180-023-01396-8. 
\bibitem{Li et al.(2018)}
Li, H., K. Yang, S. Zhao, and D. Wang. 2018. First-order random coefficients integer-valued threshold autoregressive processes. {\textit{AStA Advances in Statistical Analysis}} 102:305-331. doi:10.1007/s10182-017-0306-3.
\bibitem{Liu et al.(2015)}
Liu, Y., D. Wang, H. Zhang, and N. Shi. 2016. Bivariate zero truncated Poisson INAR(1) process. {\textit{Journal of the Korean Statistical Society}} 45(2):260-275. doi:10.1016/j.jkss.2015.11.002.
\bibitem{Monteiro et al.(2021)}
Monteiro, M., I. Pereira, and M. G. Scotto. 2021. Bivariate models for time series of counts: A comparison study between PBINAR models and dynamic factor models. {\textit{Communications in Statistics-Simulation and Computation}} 50(7):1873-1887. doi:10.1080/03610918.2019.1599015.




\bibitem[Pedeli and Karlis(2011)]{Pedeli2011}%check
Pedeli, X., and D. Karlis. 2011. A bivariate INAR(1) process with application. {\textit{Statistical Modelling}} 11(4):325-349. doi:10.1177/1471082X100110040.
\bibitem[Pedeli and Karlis(2013a)]{Pedeli2013a}%check
Pedeli, X., and D. Karlis. 2013a. On estimation of the bivariate poisson INAR process. {\textit{Communications in Statistics-Simulation and Computation}} 42(3):514-533. doi:10.1080/03610918.2011.639001.
\bibitem[Pedeli and Karlis(2013b)]{Pedeli2013b}%check
Pedeli, X., and D. Karlis. 2013b. Some properties of multivariate INAR(1) processes. {\textit{Computational Statistics and Data Analysis}} 67:213-225. doi:10.1016/j.csda.2013.05.019.
\bibitem[Popovi\'{c}(2015)]{Popović2015}%check
Popovi\'{c}, P. M. 2015. Random coefficient bivariate INAR(1) process. {\textit{Facta Universitatis, Series: Mathematics and Informatics}} 30:263-280.
\bibitem[Popovi\'{c} et al.(2016)]{Popović2016}%check
Popovi\'{c}, P. M., M. M. Risti\'{c}, and A. S. Nasti\'{c}. 2016. A geometric bivariate time series with different marginal parameters. {\textit{Statistical Papers}} 57:731-753. doi:10.1007/s00362-015-0677-z.


\bibitem[Qian et al.(2020)]{Qian2020}%check
Qian, L., Q. Li, and F. Zhu. 2020. Modelling heavy-tailedness in count time series. {\textit{Applied Mathematical Modelling}} 82:766-784. doi:10.1016/j.apm.2020.02.001.
\bibitem[Qian and Zhu(2022)]{Qian2022}%check
Qian, L., and F. Zhu. 2022. A new minification integer-valued autoregressive process driven by explanatory variables. {\textit{Australian and New Zealand Journal of Statistics}} 64(4):478-494. doi:10.1111/anzs.12379.




\bibitem[Risti\'{c} et al.(2009)]{Ristić2009}%check
Risti\'{c}, M. M., H. S. Bakouch, and A. S. Nasti\'{c}. 2009. A new geometric first-order integer-valued autoregressive (NGINAR(1)) process. {\textit{Journal of Statistical Planning and Inference}} 139(7):2218-2226. doi:10.1016/j.jspi.2008.10.007.
\bibitem[Risti\'{c} et al.(2012)]{Ristić2012}%check
Risti\'{c}, M. M., A. S. Nasti\'{c}, K. Jayakumar, and  H. S. Bakouch. 2012. A bivariate INAR(1) time series model with geometric marginals. {\textit{Applied Mathematics Letters}} 25(3):481-485. doi:10.1016/j.aml.2011.09.040.


\bibitem[Scotto et al.(2014)]{Scotto2014}
Scotto, M. G., C. H. Wei\ss, M. E. Silva, and I. Pereira. 2014. Bivariate binomial autoregressive models. {\textit{Journal of Multivariate Analysis}} 125:233-251. doi:10.1016/j.jmva.2013.12.014.
\bibitem[Su and Zhu(2021)]{Su2021}
Su, B., and F. Zhu. 2021. Comparison of BINAR(1) models with bivariate negative binomial innovations and explanatory variables. {\textit{Journal of Statistical Computation and Simulation}} 91(8):1616-1634. doi:10.1080/00949655.2020.1863965.


\bibitem{Tzougas2021}
Tzougas, G., and A. P. di Cerchiara. 2021. The multivariate mixed negative binomial regression model with an application to insurance a posteriori ratemaking. {\textit{Insurance: Mathematics and Economics}} 101: 602-625. doi:10.1016/j.insmatheco.2021.10.001.


\bibitem[Wang(1982)]{Wang1982}
Wang, Z. K. 1982. {\textit{Stochastic Process}}. Beijing: Scientific Press.


\bibitem{Wei2018}
Wei\ss, C. H. 2018. An introduction to discrete-valued time series. John Wiley \& Sons.
\bibitem{Yang2021}
Yang, K., H. Li, D. Wang, and C. Zhang. 2021. Random coefficients integer-valued threshold autoregressive processes driven by logistic regression. {\textit{AStA Adv Stat Anal}} 105:533-557. doi:10.1007/s10182-020-00379-0.
\bibitem{Yang2022}
Yang, K., X. Yu, Q. Zhang, and X. Dong. 2022. On MCMC sampling in self-exciting integer-valued threshold time series models. {\textit{Comput Stat Data Anal}} 169:107-410. doi:10.1016/j.csda.2021.107410.
\bibitem{Yang2023}
Yang, K., Y. Zhao, H. Li, and D. Wang. 2023a. On bivariate threshold Poisson integer-valued autoregressive processes. {\textit{Metrika}} 86:931-963. doi:10.1007/s00184-023-00899-0.
\bibitem{Yang2023}
Yang, K., N. Xu, H. Li, Y. Zhao, and X. Dong. 2023b. Multivariate threshold integer-valued autoregressive processes with explanatory variables. {\textit{Appl Math Model}} 124:142-166. doi:10.1016/j.apm.2023.07.030.


\bibitem[Yu et al.(2020)]{Yu2020}
Yu, M., D. Wang, K. Yang, and Y. Liu. 2020. Bivariate first-order random coefficient integer-valued autoregressive processes. {\textit{Journal of Statistical Planning and Inference}} 204:153-176. doi:10.1016/j.jspi.2019.05.004.


\bibitem[Zhang et al.(2020)]{Zhang2020}%check
Zhang, Q., D. Wang, and X. Fan. 2020. A negative binomial thinning-based bivariate INAR(1) process. {\textit{Statistica Neerlandica}} 74(4):517-537. doi:10.1111/stan.12210.
\bibitem[Zheng and Basawa(2007)]{Zheng2007}%check
Zheng, H., I. V. Basawa, and S. Datta. 2007. First-order random coefficient integer-valued autoregressive processes. {\textit{Journal of Statistical Planning and Inference}} 137(1):212-229. doi:10.1016/j.jspi.2005.12.003.
\end{thebibliography}
\end{document}